\documentclass{article}
\usepackage{mathrsfs}
\usepackage{bbm}
\usepackage{amssymb}
\usepackage{amsfonts}

\usepackage{amsmath}
\usepackage{latexsym,amsfonts,amssymb}
\newcommand\bgeq{\begin{equation}}
\newcommand\edeq{\end{equation}}
\newcommand\bgar{\begin{array}}
\newcommand\edar{\end{array}}

\title{Asymptotic stability of the stationary solution for a parabolic-hyperbolic
 free boundary problem\\ modeling tumor growth
 \footnote{This work is supported by the China
  National Natural Science Foundation under grant number 11171357.}}
\author{Shangbin Cui\footnote{E-mail:  cuisb3@yahoo.com.cn}}
\date{{\small Department of Mathematics, Sun Yat-Sen University,
   Guangzhou, Guangdong 510275,} \\ [-0.05cm]
  {\small People's Republic of China}}

\begin{document}

\maketitle

\begin{abstract}
  This paper studies asymptotic behavior of solutions of a free boundary
  problem modeling the growth of tumors with two species of cells: proliferating
  cells and quiecent cells. In previous literatures it has been proved that
  this problem has a unique stationary solution which is asymptotically stable
  in the limit case $\varepsilon=0$. In this paper we consider the more
  realistic case $0<\varepsilon\ll 1$. In this case, after suitable reduction
  the model takes the form of a coupled system of a parabolic equation and a
  hyperbolic system, so that it is more difficult than the  limit case
  $\varepsilon=0$. By using some unknown variable transform as well as the
  similarity transform technique developed in our previous work, we prove that
  the stationary solution is also asymptotically stable in the case
  $0<\varepsilon\ll 1$.


\end{abstract}
\section{Introduction}
\setcounter{equation}{0}

\hskip 2em
  It has long been observed that under a constant circumstance, a solid tumor
  will evolve into a dormant or stationary state \cite{Casey, Laird1, Laird2}.
  In a dormant state, the tumor does not change in size, while cells inside it
  are alive and keep undergoing the process of proliferation. In 1972 Greenspan
  established the first mathematical model in the form of a free boundary
  problem of a system of partial differential equations to illustrate this
  phenomenon \cite{Green1}. Since then an increasing number of tumor models in
  similar forms have been proposed by many different groups of researchers (see,
  e.g., \cite{AdmBel}, \cite{tumrev1}--\cite{ByrCha2}, \cite{tumrev2},
  \cite{Fried1}--\cite{Fried3}, \cite{PPM}--\cite{WarKin2},
  and the references cited therein). Rigorous mathematical analysis of such
  models has drawn great attention during the past twenty years, and many
  interesting results have been obtained, cf., \cite{ChenCuiF} -- \cite{CuiWei},
  \cite{FonFri1}, \cite{FriHu1} -- \cite{FriRei2}, \cite{WuCui}, and references
  cited therein.

  In this paper we study the following free boundary problem modeling the growth
  of a solid tumor with two species of cells --- proliferating cells and
  quiescent cells:
\begin{equation}
  \varepsilon c_t=\triangle{c}-F(c)
   \quad \mbox{for}\;\; x\in\Omega(t),\;\; t>0,
\end{equation}
\begin{equation}
  c=1 \quad \mbox{for}\;\; x\in\partial\Omega(t), \;\; t>0,
\end{equation}
\begin{equation}
   p_t+\nabla\cdot(\vec{v}p)=
  [K_B(c)-K_Q(c)]p+K_P(c)q \quad \mbox{for}\;\; x\in\Omega(t),\;\; t>0,
\end{equation}
\begin{equation}
   q_t+\nabla\cdot(\vec{v}q)=
  K_Q(c)p-[K_P(c)+K_D(c)]q \quad \mbox{for}\;\; x\in\Omega(t),\;\; t>0,
\end{equation}
\begin{equation}
  p+q=1 \quad \mbox{for}\;\; x\in\Omega(t),\;\; t>0,
\end{equation}
\begin{equation}
  {d\over dt}\big({\rm vol}(\Omega(t))\big)=\int_{\Omega(t)}\nabla\cdot\vec{v}(x,t)dx
  \quad {\rm for}\;\; t>0.
\end{equation}
  Here $\Omega(t)$ is the domain occupied by the tumor at time $t$, $c=c(x,t)$,
  $p=p(x,t)$ and $q=q(x,t)$ are the concentration of nutrient, the density of
  proliferating cells and the density of quiescent cells, respectively, and
  $\vec{v}=\vec{v}(x,t)$ is the velocity of tumor cell movement. Besides, $F(c)$
  is the consumption rate of nutrient by tumor cells, $K_B(c)$ is the birth rate
  of tumor cells, $K_P(c)$ and $K_Q(c)$ are transferring rates of tumor cells from
  quiescent state to proliferating state and from proliferating state to quiescent
  state, respectively, and $K_D(c)$ is the death rate of quiescent cells.

  Equation (1.1) describes diffusion of nutrient (regarded as one species) within
  the tumor, where $\varepsilon=T_{\rm diffusion}/T_{\rm growth}$ is the ratio
  of the nutrient diffusion time scale to the tumor growth (e.g., tumor doubling)
  time scale. Typically, $T_{\rm diffusion}\approx 1$min (see pp. 194--195 of
  \cite{AdmBel}), while $T_{\rm growth}\approx 1$day, so that $\varepsilon\ll 1$.
  Equation (1.2) reflects the fact that the tumor receives a constant
  nutrient supply from its boundary, and we have rescaled $c$ so that the constant
  supply of it from the tumor boundary is exactly 1. Equations (1.3) and (1.4)
  are conservation laws for the proliferating and quiescent tumor cells,
  respectively. Equation (1.5) reflects the fact that the mixture of proliferating
  and quiescent cells \underline{}in the tumor has a constant density, and we have
  rescaled $p$ and $q$ such that this constant equals to $1$. Finally, the
  equation (1.6) reflects the fact that the variance rate of the whole tumor
  volume is equal to accumulation of all local volume variance rate within
  the tumor (recall that $\nabla\cdot\vec{v}(x,t)$ is the local volume
  variance rate). We  assume that $F(c)$, $K_B(c)$, $K_P(c)$, $K_Q(c)$, $K_D(c)$
  are $C^{\infty}$ functions, and they satisfy the following conditions:
\begin{enumerate}
\item[]$(A_1)$\ \ $F'(c)>0$ for $c\in\mathbb{R}$, and $F(0)=0$.
\item[]$(A_2)$\ \ $K_B'(c)>0$ and $K_P'(c)\geq 0$ for $c\in\mathbb{R}$, and
  $K_B(0)=K_P(0)=0$.
\item[]$(A_3)$\ \ $K_D'(c)\leq 0$, $K_Q'(c)\leq 0$, $K_D(c)\geq 0$ and
  $K_Q(c)\geq 0$ for $c\in\mathbb{R}$.
\item[]$(A_4)$\ \ $K_B'(c)+K_D'(c)>0$ for $c\in\mathbb{R}$.
\item[]$(A_5)$\ \ $K_P'(c)+K_Q'(c)>0$ for $c\in\mathbb{R}$.
\end{enumerate}

  For simplicity we assume that the tumor is a strict spheroid of
  radius $R(t)$, so that
$$
  \Omega(t)=\{x\in\mathbb{R}^3: r=|x|<R(t)\}.
$$
  Accordingly, we assume that all unknown functions $c$, $p$ and $q$
  are spherically symmetric in the space variable, i.e.,
$$
  c=c(r,t), \quad p=p(r,t), \quad q=q(r,t),
$$
  and
$$
  \vec{v}=v(r,t)\frac{x}{r}.
$$
  It follows that the above model reduces into the following system of equations:
\begin{equation}
  \varepsilon c_t(r,t)=
   c_{rr}(r,t)
  +{2\over r} c_r(r,t)
  -F({c}(r,t)) \quad \mbox{for}\;\; 0<r<R(t),\;\; t>0,
\end{equation}
\begin{equation}
   c_r(0,t)=0,\quad {c}(R(t),t)=1
  \quad \mbox{for}\;\; t>0,
\end{equation}
\begin{equation}
\begin{array}{rl}
  \displaystyle p_t(r,t)+v(r,t) p_r(r,t)=&\displaystyle f({c}(r,t),p(r,t))
   \quad \mbox{for}\;\; 0<r<R(t),\;\; t>0,
\end{array}
\end{equation}
\begin{equation}
   v_r(r,t)+{2\over r}v(r,t)=g({c}(r,t),p(r,t))
  \quad \mbox{for}\;\; 0<r<R(t),\;\; t>0,
\end{equation}
\begin{equation}
  v(0,t)=0 \quad \mbox{for}\;\; t>0,
\end{equation}
\begin{equation}
  \dot{R}(t)= v(R(t),t) \quad \mbox{for}\;\; t>0,
\end{equation}
  where
$$
\begin{array}{rcl}
  f({c},p)&=&K_P(c)+\big[K_M(c)\!-\!K_N(c)\big]p-\!K_M(c)p^2,\\
  g({c},p)&=&K_M(c)p-K_D(c),
\end{array}
$$
  and
$$
  K_M(c)=K_B(c)+K_D(c), \quad
  K_N(c)=K_P(c)+K_Q(c).
$$

  Global well-posedness of this model has been proved in the literature
  \cite{CuiWei} in a more general setting (see \cite{CuiFri2} for the limit
  case $\varepsilon=0$). In \cite{CuiFri1} it is proved that under the
  assumptions $(A_1)$--$(A_5)$ the above model has a unique stationary solution.
  In \cite{ChenCuiF,Cui3} it was proved that in the limit case
  $\varepsilon=0$, this unique stationary solution is asymptotically stable.
  The purpose of this paper is to prove that the result of \cite{ChenCuiF,Cui3}
  can be extended to the more realistic case $\varepsilon\neq 0$ but small,
  namely, we shall prove that there exists $\varepsilon_0>0$ such that the
  stationary solution of the system (1.7)--(1.12) is also asymptotically stable
  in the case $0<\varepsilon<\varepsilon_0$.

  To present the precise statement of our main result, we make further
  reduction to the model as follows. Let
\begin{equation}
   \bar{{c}}(\bar{r},t)={c}(\bar{r}e^{z(t)},t), \quad \bar{p}(\bar{r},t)=
   p(\bar{r}e^{z(t)},t), \quad  \bar{v}(\bar{r},t)=
   v(\bar{r}e^{z(t)},t)e^{-z(t)}, \quad R(t)=e^{z(t)},
\end{equation}
  where $0\leq\bar{r}\leq 1$, $t\geq 0$. Then the system (1.7)--(1.12) is
  further reduced into the following system (for simplicity of the notation
  we remove all bars ``$\bar{\;\;\;}$''):
\begin{equation}
   \varepsilon e^{2z(t)}c_t( r,t)= c_{rr}( r,t)+
   \big[{2\over{r}}+\varepsilon e^{2z(t)}rv(1,t)\big]
    c_r({r},t)-e^{2z(t)}F({{c}}({r},t))
   \quad \mbox{for} \;\; 0<{r}<1,\;\; t>0,
\end{equation}
\begin{equation}
    c_r(0,t)=0,\quad {{c}}(1,t)=1
   \quad  \mbox{for}\;\; t>0,
\end{equation}
\begin{eqnarray}
   \displaystyle p_t({r},t)+{w}({r},t) p_r({r},t)
   =&&f({c}({r},t),{p}({r},t))
   \quad \mbox{for} \;\; 0<{r}<1,\;\; t>0,
\end{eqnarray}
\begin{equation}
    v_r({r},t)
   +{2\over{r}}{v}({r},t)
   =g({c}({r},t),{p}({r},t))
   \quad \mbox{for} \;\; 0<{r}<1,\;\; t>0,
\end{equation}
\begin{equation}
   {v}(0,t)=0 \quad  \mbox{for}\;\; t>0,
\end{equation}
\begin{equation}
   \dot{z}(t)={v}(1,t)  \quad \mbox{for} \;\; t>0,
\end{equation}
  where
\begin{equation}
   {w}({r},t)={v}({r},t)-{r}{v}(1,t).
\end{equation}

  Let $(c^*(r),p^*(r),v^*(r),z^*)$ be the stationary solution of the system
  (1.14)--(1.20) ensured by \cite{CuiFri1}. The main result of this paper is
  the following:
\medskip

  {\bf Theorem 1.1}\ \ {\em Let the assumptions $(A_1)$--$(A_5)$ be satisfied.
  There exist positive constants $\varepsilon_0$ and $\mu^*$ such that for any
  $0<\varepsilon<\varepsilon_0$ and $0<\mu<\mu^*$, there exist corresponding
  positive constants $\delta$ and $C$ such that for any time-dependent solution
  $(c(r,t),p(r,t),v(r,t),z(t))$ of the system $(1.14)$--$(1.20)$, if the initial
  data $(c_0(r),p_0(r),z_0)=(c(r,0),p(r,0),z(0))$ satisfy
\begin{equation}
  c_0(r),\;c_0'(r)\in C[0,1], \quad c_0'(0)=0, \quad c_0(1)=1, \quad
  0\leq c_0(r)\leq 1  \quad \mbox{for}\;\;0\leq r\leq 1,
\end{equation}
\begin{equation}
  p_0(r),\;r(1-r)p_0'(r)\in C[0,1], \quad p_0(1)=1,  \quad 0\leq p_0(r)\leq 1
   \quad \mbox{for}\;\;0\leq r\leq 1,
\end{equation}
\begin{equation}
  \max_{0\leq r\leq 1}|c_0(r)-c^*(r)|<\delta,
  \quad
  \sup_{0\leq r\leq 1}|c_0'(r)-c^{*'}(r)|<\delta,
\end{equation}
\begin{equation}
  \max_{0\leq r\leq 1}|p_0(r)-p^*(r)|<\delta,
  \quad
  \sup_{0\leq r\leq 1}r(1\!-\!r)|p_0'(r)-p^{*'}(r)|<\delta,
  \quad \mbox{and} \quad
  |z_0-z^*|<\delta,
\end{equation}
  then for all $t\geq 0$ we have
\begin{equation}
  \max_{0\leq r\leq 1}[|c(r,t)-c^*(r)|+
  \max_{0\leq r\leq 1}|c_r(r,t)-c^{*'}(r)|+
  \max_{0\leq r\leq 1}|c_t(r,t)|]<C\delta e^{-\mu t},
\end{equation}
\begin{equation}
  \max_{0\leq r\leq 1}|p(r,t)-p^*(r)|+
  \max_{0\leq r\leq 1}r(1\!-\!r)|p_r(r,t)-p^{*'}(r)|+
  \max_{0\leq r\leq 1}|p_t(r,t)|]<C\delta e^{-\mu t},
\end{equation}
  and
\begin{equation}
  |z(t)-z_*|+|\dot{z}(t)|<C\delta e^{-\mu t}.
\end{equation}
}

  {\em Remark}:\ \ Readers who are familiar with the literatures \cite{FriRei1},
  \cite{CuiFri0} and \cite{Cui2} might wish to use the iteration technique
  used there to prove the above theorem. The author is not successful in trying
  to do so. The obstacle lies in making estimates for $|\bar{r}-S(\bar{r},t,s)|$
  and related quantities, where $S(\bar{r},t,s)$ denotes the similarity transform
  introduced in \cite{Cui3} (see (4.16) in Section 4 below). In this paper we
  use a different approach; see Section 2 for illustration.

  The structure of the rest part is as follows: In the next section we make
  further reduction to the system (1.14)--(1.20) and illustrate the main idea of
  the proof of the above theorem. In Sections 3--5 we make preparations for
  the proof of Theorem 1.1, which is given in the last section.

\section{Reduction of the problem}
\setcounter{equation}{0}

\hskip 2em
  We define a function $ m$ in $[0,1]\times\mathbb{R}$ as follows: For
  any $z\in\mathbb{R}$, the function $c(r)= m(r;z)$ $(0\leq r\leq 1)$
  is the unique solution of the following boundary value problem:
\begin{equation}
\left\{
\begin{array}{l}
   \displaystyle {c}''(r)+{2\over r}{c}'(r)=e^{2z}F({c}(r))
   \quad \mbox{for} \;\; 0<r<1,
\\ [0.3cm]
   \displaystyle {c}'(0)=0,\quad {c}(1)=1.
\end{array}
\right.
\end{equation}
  Since $F\in C^{\infty}(\mathbb{R})$, we have $ m\in C^{\infty}([0,1]
  \times\mathbb{R})$. For the solution $(c(r,t),p(r,t),v(r,t),z(t))$ of the problem
  (1.14)--(1.19), we let
$$
  \eta(r,t)=c(r,t)- m(r;z(t)).
\eqno{(2.2)}
$$
  It can be easily seen that $\eta$ satisfies the following equations:
$$
\left\{
\begin{array}{l}
   \displaystyle\varepsilon e^{2z(t)}\eta_t(r,t)=\eta_{rr}(r,t)
  +\big[{2\over r}+\varepsilon e^{2z(t)}rv(1,t)\big]\eta_r(r,t)
  -e^{2z(t)}a(r;\eta(r,t),z(t))\eta(r,t)
\\ [0.1cm]
   \quad\quad\quad\quad\quad
   -\varepsilon e^{2z(t)}v(1,t)\big[ m_z(r;z(t))
  -r m_r(r;z(t))\big]\quad \mbox{for}\;\; 0<r<1,\;\; t>0,
\\ [0.1cm]
   \displaystyle\eta_r(0,t)=0,\quad \eta(1,t)=0
  \quad \mbox{for}\;\; t>0,
\end{array}
\right.
\eqno{(2.3)}
$$
  where $a$ is a function in $[0,1]\times\mathbb{R}\times\mathbb{R}$
  defined as follows: For $r\in [0,1]$ and $y,z\in\mathbb{R}$,
$$
  a(r;y,z)=\int_0^1F'( m(r;z)+\theta y)d\theta.
\eqno{(2.4)}
$$
  Conversely, if $\eta(r,t)$ satisfies (2.3) and $z(t)$ satisfies (1.19), then
  $c(r,t)= m(r;z(t))+\eta(r,t)$ satisfies (1.14) and (1.15). Hence,
   under the transformation of unknown variables
$$
  (c(r,t),p(r,t),v(r,t),z(t))\mapsto(\eta(r,t),p(r,t),v(r,t),z(t))
\eqno{(2.5)}
$$
  given by (2.2), the system (1.14)--(1.20) is equivalent to the system
  (1.16)--(1.20) (with $c(r,t)$ replaced by $ m(r;z(t))+\eta(r,t)$)
  coupled with (2.3).
\medskip

  We introduce three maps $\mathbbm{v},\mathbbm{w},\mathbbm{f}: C[0,1]\times
  C[0,1]\times\mathbb{R}\rightarrow C[0,1]$ and a functional $\mathbbm{g}$ in
  $C[0,1]\times C[0,1]\times\mathbb{R}$ as follows: For $\eta,p\in C[0,1]$ and
  $z\in\mathbb{R}$,
$$
\begin{array}{rl}
  \mathbbm{v}(r;\eta,p,z)=&\displaystyle\left\{
\begin{array}{l}
   \displaystyle\frac{1}{r^2}\int_0^rg( m(\rho,z)+\eta(\rho),p(\rho))\rho^2d\rho
   \quad \mbox{for} \;\; 0<r\leq 1,
\\ [0.2cm]
   0,\quad \mbox{for} \;\; r=0;
\end{array}
\right.
\\ [0.2cm]
  \mathbbm{w}(r;\eta,p,z)=&\mathbbm{v}(r;\eta,p,z)-r\mathbbm{v}(1;\eta,p,z),
\\ [0.2cm]
  \mathbbm{f}(r;\eta,p,z)=&f( m(r,z)+\eta(r),p(r)),
\\
  \mathbbm{g}(\eta,p,z)=&\mathbbm{v}(1;\eta,p,z)
  =\displaystyle\int_0^1g(m(r,z)+\eta(r),p(r))r^2dr.
\end{array}
$$
  It can be easily seen that with $\eta$ as in (2.2), we have
$$
\begin{array}{rcl}
   v(r,t)&=&\mathbbm{v}(r;\eta(\cdot,t),p(\cdot,t),z(t)),
\\ [0.1cm]
   w(r,t)&=&\mathbbm{w}(r;\eta(\cdot,t),p(\cdot,t),z(t)),
\\ [0.1cm]
   f(c(r,t),p(r,t))&=&\mathbbm{f}(r;\eta(\cdot,t),p(\cdot,t),z(t)),
\\ [0.1cm]
   v(1,t)&=&\mathbbm{g}(\eta(\cdot,t),p(\cdot,t),z(t)).
\end{array}
$$
   Hence, under the transformation of unknown variables given by (2.2) and (2.5),
   the  system (1.14)--(1.20) reduces into the following system:
$$
\left\{
\begin{array}{l}
   \displaystyle\varepsilon e^{2z(t)}\eta_t(r,t)=\eta_{rr}(r,t)
  +\big[{2\over r}+\varepsilon e^{2z(t)}r\mathbbm{g}(\eta,p,z)\big]\eta_r(r,t)
  -e^{2z(t)}a(r;\eta,z)\eta
\\
   \quad\quad\quad\quad\quad
   -\varepsilon e^{2z(t)}\mathbbm{g}(\eta,p,z)\big[ m_z(r;z)
  -r m_r(r;z)\big]\quad \mbox{for}\;\; 0<r<1,\;\; t>0,
\\
   \displaystyle\eta_r(0,t)=0,\quad \eta(1,t)=0
  \quad \mbox{for}\;\; t>0,
\\
   \displaystyle p_t
   +\mathbbm{w}(r;\eta,p,z) p_r
   =\mathbbm{f}(r;\eta,p,z) \quad \mbox{for} \;\; 0<r<1,\;\; t>0,
\\
   \displaystyle\dot{z}=\mathbbm{g}(\eta,p,z) \quad \mbox{for} \;\; t>0.
\end{array}
\right.
\eqno{(2.6)}
$$

  We shall treat the above system in the following way: Let
$$
\begin{array}{c}
  Y=\displaystyle\Big\{\eta=\eta(r,t)\in C([0,1]\times [0,\infty)):\,
  \eta(r,t)\;\mbox{is differentiable in $r$},\;
  \eta_r\in C([0,1]\times [0,\infty)),
\\ [0.2cm]
  \displaystyle\quad\;\;
  \eta_r(0,t)=\eta(1,t)=0,\;\forall t\geq 0,\;\mbox{and}\;
  \|\eta\|_Y\stackrel{\mathrm{def}}{=}
  \sup_{0\leq r\leq 1\atop t\geq 0}e^{\mu t}
  |\eta(r,t)|+\sup_{0\leq r\leq 1\atop t\geq 0}e^{\mu t}
  \Big|\eta_r(r,t)\Big|<\infty\Big\},
\end{array}
$$
  where $\mu$ is a positive constant to be specified later. It is clear that
  $(Y,\|\cdot\|_Y)$ is a Banach space. Let $\delta$ and $\delta'$ be two positive
  constants to be specified later. We assume that the initial value
  $(c_0(r),p_0(r),z_0)$ of $(c(r,t),p(r,t),z(t))$ satisfies the conditions in
  Theorem 1.1. Given $\eta\in Y$ satisfying the conditions
$$
  \|\eta\|_Y\leq\delta' \quad \mbox{and} \quad
  \eta(r,0)=\eta_0(r)\stackrel{\mathrm{def}}{=}c_0(r)- m(r,z_0),
$$
  we first solve the initial value problem
$$
\left\{
\begin{array}{l}
   \displaystyle p_t
   +\mathbbm{w}(r;\eta,p,z) p_r
   =\mathbbm{f}(r;\eta,p,z)
   \quad \mbox{for} \;\; 0<r<1\;\; \mbox{and}\;\; t>0,
\\ [0.2cm]
   \displaystyle\dot{z}=\mathbbm{g}(\eta,p,z)
   \quad \mbox{for} \;\; t>0,
\\ [0.1cm]
   p(r,0)=p_0(r) \quad \mbox{for} \;\; 0\leq r\leq 1, \quad \mbox{and}
    \quad z(0)=z_0,
\end{array}
\right.
\eqno{(2.7)}
$$
  and next solve the following initial-boundary value problem:
$$
\left\{
\begin{array}{l}
   \displaystyle\varepsilon e^{2z(t)}\tilde{\eta}_t=
  \tilde{\eta}_{rr}
  +\big[{2\over r}+\varepsilon e^{2z(t)}r\mathbbm{g}(\eta,p,z)\big]
  \tilde{\eta}_r
  -e^{2z(t)}a(r;\eta,z)\tilde{\eta}
\\ [0.2cm]
   \quad\quad\quad\quad\quad
   -\varepsilon e^{2z(t)}\mathbbm{g}(\eta,p,z)
   \big[m_z(r;z)-r m_r(r;z)\big]
   \quad \mbox{for}\;\; 0<r<1,\;\; t>0,
\\ [0.1cm]
   \displaystyle\tilde{\eta}_r(0,t)=0,
   \quad \tilde{\eta}(1,t)=0 \quad \mbox{for}\;\; t>0,
\\ [0.1cm]
   \tilde{\eta}(r,0)=\eta_0(r) \quad \mbox{for} \;\; 0\leq r\leq 1.
\end{array}
\right.
\eqno{(2.8)}
$$
  We then obtain a map $\eta\mapsto\tilde{\eta}$ from the closed ball
  $\overline{B}_{\delta'}(0)$ in $Y$ into $Y$. We shall prove that there exists
  $\varepsilon_0>0$ such that for any $0<\varepsilon\leq\varepsilon_0$, if
  $\delta$ and $\delta'$ are sufficiently small then this is a self-mapping in
  $\overline{B}_{\delta'}(0)$. Moreover, we shall use the Schauder fixed
  point theorem to prove it has a fixed point in this ball. The desired
  assertion then follows.

\section{A preliminary lemma}
\setcounter{equation}{0}

\hskip 2em
  {\bf Lemma 3.1}\ \ {\em Assume that $F'(c)>0$ for $c\in\mathbb{R}$ and
  $F(0)=0$. For all $0\leq r\leq 1$ and $z\in\mathbb{R}$ we have:

  $(1)$ $0<m(r;z)\leq 1$; $\qquad\qquad\;\;$
  $(2)$ $m_r(r;z)\geq 0$; $\qquad\qquad\;\;$
  $(3)$ $m_z(r;z)\leq 0$;

  $(4)$ $m_r(r;z)\leq\displaystyle\frac{1}{3}rF(1)e^{2z}$; $\qquad$
  $(5)$ $m_z(r;z)\geq -\displaystyle\frac{1}{3}(1-r^2)F(1)e^{2z}$.

  $(6)$ $\displaystyle-\frac{2}{3}F(1)e^{2z}\leq m_{rr}(r;z)\leq\displaystyle F(1)e^{2z}$;

  $(7)$ $m_z(r;z)\geq -\displaystyle\frac{1}{3}(1-r^2)F(1)e^{2z}$.
}
\medskip

  {\em Proof}:\ \ The assertions (1)$\sim$(3) are well-known. Let $c(r)=
   m(r;z)$. From the first equation in (2.1) we have
$$
  \big(c'(r)r^2\big)'=e^{2z}F(c(r))r^2.
$$
  Integrating both sides from $0$ to $r$ and using the assertion (1) we obtain
\begin{equation}
  c'(r)=\displaystyle\frac{e^{2z}}{r^2}\int_0^rF(c(\rho))\rho^2d\rho
  \leq\frac{1}{3}rF(1)e^{2z}.
\end{equation}
  This proves the assertion (4). Next, differentiating all three equations in
  (2.1) in $z$ we see that the function $c_z(r)= m_z(r;z)$ is the
  solution of the following boundary value problem:
$$
\left\{
\begin{array}{l}
   \displaystyle {c}_z''(r)+{2\over r}{c}_z'(r)=e^{2z}F'({c}(r)){c}_z(r)
   +2e^{2z}F({c}(r)) \quad \mbox{for} \;\; 0<r<1,
\\ [0.3cm]
   \displaystyle {c}_z'(0)=0,\quad {c}_z(1)=0.
\end{array}
\right.
$$
  Since $F'(c)>0$, ${c}_z(r)\leq 0$ and ${c}(r)\leq 1$, we see that
$$
  \big(c_z'(r)r^2\big)'=e^{2z}F'({c}(r)){c}_z(r)r^2
   +2e^{2z}F({c}(r))r^2\leq 2e^{2z}F(1)r^2.
$$
  It follows that
$$
  c_z'(r)\leq\displaystyle\frac{2e^{2z}}{r^2}F(1)\int_0^r\rho^2d\rho
  \leq\frac{2}{3}rF(1)e^{2z}.
$$
  Hence
$$
  -c_z(r)\leq\displaystyle\frac{2}{3}F(1)e^{2z}\int_r^1\rho d\rho
  =\frac{1}{3}(1-r^2)F(1)e^{2z}.
$$
  This proves the assertion (5) and completes the proof of Lemma 3.1.
  $\quad\Box$

\section{Decay estimates of the solution of (2.7)}
\setcounter{equation}{0}

\hskip 2em

  We denote by $C_\vee^1[0,1]$ the following function space:
$$
  C_\vee^1[0,1]=\{\phi\in C[0,1]\cap C^1(0,1):r(1-r)\phi'(r)\in C[0,1]\},
$$
  with norm
$$
  \|\phi\|_{C_\vee^1[0,1]}=\max_{0\leq r\leq 1}|\phi(r)|+
  \sup_{0<r<1}|r(1-r)\phi'(r)|
  \quad \mbox{for} \;\; \phi\in C_\vee^1[0,1].
$$
  It is clear that this is a Banach space. Next we denote
$$
   X=C[0,1]\times\mathbb{R} \quad \mbox{and} \quad
   X_0=C_\vee^1[0,1]\times\mathbb{R};
$$
  they are both Banach spaces. Given $\eta\in Y$, we introduce a map
  ${\mathbb F}: [0,\infty)\times X_0\to X$ as follows:
  For $t\geq 0$ and $U=(p,z)\in X_0$,
$$
  {\mathbb F}(t,U)=\big(-\mathbbm{w}(\cdot;\eta(\cdot,t),p,z)p'
   +\mathbbm{f}(\cdot;\eta(\cdot,t),p,z),\mathbbm{g}(\eta(\cdot,t),p,z)\big).
$$
  Then the system (2.7) can be rewritten as the following differential equation
  in the Banach space $X$ (regarding $X_0$ as a subspace of
  $X$ and ${\mathbb F}(\cdot,t)$ as an unbounded nonlinear operator in
  $X$ with domain $X_0$):
$$
  \dot{U}={\mathbb F}(t,U) \quad \mbox{for}\;\; t>0.
\eqno{(4.1)}
$$

  It is clear that $[t\mapsto {\mathbb F}(t,\cdot)]\in C([0,\infty),
  C^{\infty}(X_0,X))$, and
$$
  {\mathbb F}(t,U)={\mathbb A}_0(t,U)U+{\mathbb F}_0(t,U)
  \quad \mbox{for}\;\; U\in X_0, \;\; t\geq 0,
$$
  where
$$
\begin{array}{rl}
  {\mathbb A}_0(t,U)V=&\big(-\mathbbm{w}(\cdot;\eta(\cdot,t),p,z)q',0\big)
   \quad \mbox{for}\;\;U=(p,z)\in X, \;\; V=(q,y)\in X_0, \;\; t\geq 0,
\\
  {\mathbb F}_0(t,U)=&\big(\mathbbm{f}(\cdot;\eta(\cdot,t),p,z),
  \mathbbm{g}(\eta(\cdot,t),p,z)\big)
   \quad \mbox{for}\;\;U=(p,z)\in X, \;\; t\geq 0.
\end{array}
$$
  It is also clear that $[t\mapsto {\mathbb A}_0(t,\cdot)]\in C([0,\infty),
  C^{\infty}(X,\mathscr{L}(X_0,X)))$ and $[t\mapsto {\mathbb F}(t,\cdot)]\in
  C([0,\infty),\\ C^{\infty}(X,X))$. Let $U^*=(p^*,z^*)$ and $V=U-U^*$. Then
$$
\begin{array}{rl}
  {\mathbb F}(t,U)=&{\mathbb A}_0(t,U)U+{\mathbb F}_0(t,U)
\\
  =&{\mathbb A}_0(t,U^*+V)(U^*+V)+{\mathbb F}_0(t,U^*+V)
\\
  =&{\mathbb A}_0(t,U^*+V)V+{\mathbb A}_0(t,U^*+V)U^*+{\mathbb F}_0(t,U^*+V)
\\
  =&\{{\mathbb A}_0(t,U^*+V)V+[\partial_U{\mathbb A}_0^*(U^*)V]U^*
  +\partial_U{\mathbb F}_0^*(U^*)V\}
\\
  &+\{[\partial_U{\mathbb A}_0(t,U^*)V-\partial_U{\mathbb A}_0^*(U^*)V]U^*
  +[\partial_U{\mathbb F}_0(t,U^*)V-\partial_U{\mathbb F}_0^*(U^*)V]\}
\\
   &+\{[{\mathbb A}_0(t,U^*+V)-{\mathbb A}_0(t,U^*)-
   \partial_U{\mathbb A}_0(t,U^*)V]U^*+[{\mathbb F}_0(t,U^*+V)
\\
   &\quad -{\mathbb F}_0(t,U^*)-\partial_U{\mathbb F}_0(t,U^*)V]\}
  +\{{\mathbb A}_0(t,U^*)U^*+{\mathbb F}_0(t,U^*)\}
\\
  =&{\mathbb A}(t,V)V+{\mathbb B}(t)V+{\mathbb G}(t,V)+G(t),
\end{array}
$$
  where $\partial_U{\mathbb A}_0(t,U)$, $\partial_U{\mathbb F}_0(t,U)$,
  $\partial_U{\mathbb A}_0^*(U)$, $\partial_U{\mathbb F}_0^*(U)$ denote
  the Fr\'{e}chet derivatives of ${\mathbb A}_0(t,U)$, ${\mathbb F}_0(t,U)$,
  ${\mathbb A}_0^*(U)$, ${\mathbb F}_0^*(U)$ in the variable $U$,
  respectively, where
$$
\begin{array}{rl}
  {\mathbb A}_0^*(U)V=&\big(-\mathbbm{w}(\cdot;0,p,z)q',0\big)
  =\displaystyle\lim_{t\to\infty}{\mathbb A}_0(t,U)V
   \quad \mbox{for}\;\;U=(p,z)\in X, \;\; V=(q,y)\in X_0,
\\
  {\mathbb F}_0^*(U)=&\big(\mathbbm{f}(\cdot;0,p,z),
  \mathbbm{g}(0,p,z)\big)=\displaystyle\lim_{t\to\infty}{\mathbb F}_0(t,U)
   \quad \mbox{for}\;\;U=(p,z)\in X,
\end{array}
$$
  and for $V\in X$, $W\in X_0$ and $t\geq 0$,
$$
\begin{array}{rl}
  {\mathbb A}(t,V)W=&{\mathbb A}_0(t,U^*+V)W+[\partial_U{\mathbb A}_0(t,U^*)W]
  U^*+\partial_U{\mathbb F}_0(t,U^*)W,
\\
  {\mathbb B}(t)V=&[\partial_U{\mathbb A}_0(t,U^*)V-\partial_U{\mathbb A}_0^*(U^*)V]U^*
  +[\partial_U{\mathbb F}_0(t,U^*)V-\partial_U{\mathbb F}_0^*(U^*)V],
\\
  {\mathbb G}(t,V)=&[{\mathbb A}_0(t,U^*+V)-{\mathbb A}_0(t,U^*)-
  \partial_U{\mathbb A}_0(t,U^*)V]U^*+[{\mathbb F}_0(t,U^*+V)
\\
  &-{\mathbb F}_0(t,U^*)-\partial_U{\mathbb F}_0(t,U^*)V],
\\
  G(t)=&{\mathbb A}_0(t,U^*)U^*+{\mathbb F}_0(t,U^*)={\mathbb F}(t,U^*).
\end{array}
$$
  Hence (4.1) can be rewritten as follows:
$$
  \dot{V}={\mathbb A}(t,V)V+{\mathbb B}(t)V+{\mathbb G}(t,V)+G(t)
  \quad \mbox{for}\;\; t>0.
\eqno{(4.2)}
$$

  Note that the above equation is quasi-linear, and
$$
  {\mathbb A}(t,V)={\mathbb A}_0(t,U^*+V)+{\mathbb B}^*,
$$
  where
$$
  {\mathbb B}^*W=[\partial_U{\mathbb A}_0^*(U^*)W]U^*
  +\partial_U{\mathbb F}_0^*(U^*)W \quad \mbox{for}\;\; W\in X.
$$
  It is easy to see that ${\mathbb B}^*\in\mathscr{L}(X)\cap\mathscr{L}(X_0)$
  (cf. Corollary 3.2 of \cite{Cui3}). Moreover, using the mean value theorem we
  can easily prove that ${\mathbb B}\in C([0,\infty),\mathscr{L}(X)\cap
  \mathscr{L}(X_0))$, and there exists constant $C>0$ independent of $\eta$
  such that
$$
  \|{\mathbb B}(t)\|_{\mathscr{L}(X)}+\|{\mathbb B}(t)\|_{\mathscr{L}(X_0)}
  \leq C[\max_{0\leq r\leq 1}|\eta(r,t)|+\max_{0\leq r\leq 1}|\eta_r(r,t)|]
  \leq C\delta e^{-\mu t}
\eqno{(4.3)}
$$
  for all $t\geq 0$. Besides, since $\eta$ and $\eta_r$ are bounded functions,
  It is also easy to see that there is a positive constant $C$ independent of
  $\eta$ such that
$$
  \|{\mathbb G}(t,V)\|_X\leq C\|V\|_X^2,
\eqno{(4.4)}
$$
$$
  \|{\mathbb G}(t,V)\|_{X_0}\leq C\|V\|_{X_0}^2,
\eqno{(4.5)}
$$
$$
  \|{\mathbb G}(t,V_1)-{\mathbb G}(t,V_2)\|_X\leq
  C(\|V_1\|_X+\|V_2\|_X)\|V_1-V_2\|_X
\eqno{(4.6)}
$$
$$
  \|{\mathbb G}(t,V_1)-{\mathbb G}(t,V_2)\|_{X_0}\leq
  C(\|V_1\|_{X_0}+\|V_2\|_{X_0})\|V_1-V_2\|_{X_0}
\eqno{(4.7)}
$$
  for all $t\geq 0$ and small $\|V\|_X,\|V_1\|_{X_0},\|V_2\|_{X_0}$ (cf.
  Corollary 3.3 of \cite{Cui3}). As for the last term on the right-hand side
  of (4.2), we have:
\medskip

  {\bf Lemma 4.1}\ \ {\em Let $\eta\in Y$ and $\|\eta\|_Y\leq\delta$. Then
  $G\in C([0,\infty),X_0)$ and there exists a positive constant $C$ indepedent
  of $\delta$ such that}
$$
  \|G(t)\|_{X_0}\leq C\delta e^{-\mu t} \quad \mbox{for}\;\; t\geq 0.
\eqno{(4.8)}
$$

  {\em Proof}:\ \ The assertion that $G\in C([0,\infty),X_0)$ follows immediately
  from Lemma 3.1 of \cite{Cui3}. To prove (4.8) we note that
$$
  G(t)={\mathbb F}(t,U^*)=\big(-\mathbbm{w}(\cdot;\eta(\cdot,t),p^*,z^*)p^{*'}
   +\mathbbm{f}(\cdot;\eta(\cdot,t),p^*,z^*),\mathbbm{g}(\eta(\cdot,t),p^*,z^*)\big).
$$
  Since $(c^*,p^*,v^*,z^s)$ is stationary solution of (1.14)--(1.20), we have
$$
\begin{array}{c}
  c^*(r)=m(r,z^*), \quad -v^*(r)p^{*'}(r)+f(c^*(r),p^*(r))=0, \quad
  \mbox{and}
\\
  v^*(1)=\displaystyle\int_0^1g(c^*(s),p^*(s))s^2ds=0.
\end{array}
$$
  Moreover, since $v^*(1)=0$, we have
$$
\begin{array}{rl}
  |\mathbbm{w}(r;\eta(\cdot,t),p^*,z^*)-v^*(r)|\leq &
  |\mathbbm{v}(r;\eta(\cdot,t),p^*,z^*)-v^*(r)|
  +r|\mathbbm{v}(1;\eta(\cdot,t),p^*,z^*)-v^*(1)|
\\ [0.3cm]
  \leq &\displaystyle\frac{1}{r^2}\int_0^r|g(c^*(\rho)+\eta(\rho,t),p^*(\rho))
  -g(c^*(\rho),p^*(\rho))|\rho^2d\rho
\\ [0.3cm]
  &\displaystyle+r\int_0^1|g(c^*(\rho)+\eta(\rho,t),p^*(\rho))
  -g(c^*(\rho),p^*(\rho))|\rho^2d\rho
\\
  \leq &\displaystyle Cr\max_{0\leq r\leq 1}|\eta(r,t)|
  \quad \mbox{for}\;\; 0\leq r\leq 1,\;\;  t\geq 0.
\end{array}
$$
  Hence
$$
\begin{array}{rl}
  \|G(t)\|_X=&\displaystyle\max_{0\leq r\leq 1}
  |\mathbbm{w}(r;\eta(\cdot,t),p^*,z^*)p^{*'}(r)-\mathbbm{f}(\cdot;\eta(\cdot,t),p^*,z^*)|
  +|\mathbbm{g}(\eta(\cdot,t),p^*,z^*)|
\\
  \leq &\displaystyle\max_{0\leq r\leq 1}
  r^{-1}|\mathbbm{w}(r;\eta(\cdot,t),p^*,z^*)-v^*(r)|\cdot\max_{0\leq r\leq 1} rp^{*'}(r)
\\
  &\displaystyle+\max_{0\leq r\leq 1}|f(c^*(r)+\eta(r,t),p^*(r))-f(c^*(r),p^*(r))|
  +|\mathbbm{v}(1;\eta(\cdot,t),p^*,z^*)-v^*(1)|
\\
  \leq &\displaystyle C\max_{0\leq r\leq 1}|\eta(r,t)|
  \leq C\delta e^{-\mu t}
  \quad \mbox{for}\;\; t\geq 0.
\end{array}
$$
  Similarly, by using the above-mentioned equations for $(c^*,p^*,v^*,z^s)$ and
  the fact that $rp^{*'}(r),\\ r^2p^{*''}(r)\in C[0,1]$ (cf. Lemma 3.1 of
  \cite{Cui3}) we can also prove that
$$
\begin{array}{rl}
  &\displaystyle\max_{0\leq r\leq 1}
  |r(1-r)\frac{\partial}{\partial r}[\mathbbm{w}(r;\eta(\cdot,t),p^*,z^*)p^{*'}(r)
  -\mathbbm{f}(\cdot;\eta(\cdot,t),p^*,z^*)]|
\\
  \leq &\displaystyle C\max_{0\leq r\leq 1}|\eta(r,t)|
  + C\max_{0\leq r\leq 1}|\eta_r(r,t)|
  \leq C\delta e^{-\mu t} \quad \mbox{for}\;\; t\geq 0.
\end{array}
$$
  Hence (4.8) holds. $\quad\Box$
\medskip

  Let $V\in C([0,\infty),X)$ be given. By using some similar arguments as in the
  proof of Lemma 4.2 of \cite{Cui3}, we can easily prove that the family of
  unbounded linear operators $\{{\mathbb A}(t,V):\; t\geq 0\}$ in $X$ is a stable
  family of infinitesimal generators of $C_0$ semigroups in $X$, and
  $\{\tilde{\mathbb A}(t,V):\; t\geq 0\}$, their parts in $X_0$, is a stable
  family of infinitesimal generators of $C_0$ semigroups in $X_0$. By Theorem 3.1
  in Chapter 5 of \cite{Pazy}, it follows that ${\mathbb A}(t,V)\; (t\geq 0)$
  generates an evolution system $\mathbb{U}(t,s;V)$ ($t\geq s\geq 0$). Moreover,
  by using a similar argument as in the proof of Lemma 4.3 of \cite{Cui3}, we can
  easily prove that for any $U_0\in X_0$ and $s\geq 0$, the initial value problem
$$
  \dot{U}(t)={\mathbb A}(t,V)U(t)\quad (\mbox{for}\;\; t>s)
  \quad \mbox{and} \quad  U(s)=U_0
\eqno{(4.9)}
$$
  has a unique solution $U=U(t;s,U_0)\in C([0,\infty),X_0)\bigcap C^1([0,\infty),X)$.
  By Theorem 4.2 in Chapter 5 of \cite{Pazy}, it follows that $U(t;s,U_0)=
  \mathbb{U}(t,s;V)U_0$ for all $t\geq s\geq 0$, and, furthermore, for any $F\in
  C([0,\infty),X_0)$ the expression
$$
  U(t)=\mathbb{U}(t,0;V)U_0+\int_0^t\mathbb{U}(t,s;V)F(s)ds
  \quad \mbox{for}\;\; t\geq0.
$$
  gives a unique solution $U\in C([0,\infty),X_0)\bigcap C^1([0,\infty),X)$ of the
  initial value problem
$$
  \dot{U}(t)={\mathbb A}(t,V)U(t)+F(t)\quad (\mbox{for}\;\; t>0)
  \quad \mbox{and} \quad  U(0)=U_0.
$$
  Hence, given $V_0\in X_0$, the unique solution of the equation (4.2) subject
  to the initial condition $V(0)=V_0$ is also a solution of the integral equation
$$
  V(t)=\mathbb{U}(t,0;V)V_0+\int_0^t\mathbb{U}(t,s;V)[\mathbb{B}(s)V
  +\mathbb{G}(s,V)ds+G(s)]ds
\eqno{(4.10)}
$$
  and vice versa, i.e., the equation (4.2) subject to the initial condition $V(0)
  =V_0$ is equivalent to the above integral equation.
\medskip

  {\bf Lemma 4.2}\ \ {\em Let $\eta\in C([0,1]\times [0,\infty))$ and $V=V(t)\in
  C([0,\infty),X)$ be given such that $\eta_r\in C([0,1]\times [0,\infty))$,
$$
  \sup_{0\leq r\leq 1}|\eta(r,t)|+\sup_{0\leq r\leq 1}|\eta_r(r,t)|
  \leq\delta e^{-\mu t} \quad \mbox{and} \quad
  \|V(t)\|_{X}\leq C_0\delta e^{-\mu t} \quad \mbox{for} \;\; t\geq 0,
\eqno{(4.11)}
$$
  where $C_0$, $\delta$ and $\mu$ are given positive constants. There exists a
  positive constant $\mu^*$ determined by the operator $\mathbb{A}^*(0)=
  \displaystyle\lim_{t\to\infty}\mathbb{A}(t,U^*)$ $($therefore it is independent
  of the constants $C_0$, $\delta$ and $\mu)$, such that for any $0<\mu'<\mu^*$
  there exists corresponding constant $\delta_0>0$ $($depending on $C_0$, $\mu$
  and $\mu')$ such that if $0<\delta\leq\delta_0$ then the following estimates
  hold:
$$
  \|{\mathbb U}(t,s, V)\|_{L(X)}\leq C_1 e^{-\mu'(t-s)} \quad
  \mbox{for} \;\; t\geq s\geq 0,
\eqno{(4.12)}
$$
$$
  \|{\mathbb U}(t,s, V)\|_{L(X_0)}\leq C_2 e^{-\mu'(t-s)} \quad
  \mbox{for} \;\; t\geq s\geq 0,
\eqno{(4.13)}
$$
  where $C_1$ and $C_2$ are positive constants depending only on $\mu'$ and not
  depending on $\mu$ and $C_0$.}
\medskip

  {\em Proof}:\ \ Let $V(t)=(\varphi(\cdot,t),\zeta(t))$ and set
$$
   p(r,t)=p^*(r)+\varphi(r,t), \quad  z(t)=z^*+\zeta(t).
$$
  We denote
$$
\begin{array}{rl}
  w(r,t)=&\mathbbm{w}(r;\eta(\cdot,t),p(\cdot,t),z(t))
\\ [0.2cm]
  =&\displaystyle\frac{1}{r^2}\int_0^rg(m(\rho,z(t))+\eta(\rho,t),
   p(\rho,t))\rho^2d\rho-r\int_0^1g(m(\rho,z(t))+\eta(\rho,t),
   p(\rho,t))\rho^2d\rho
\end{array}
$$
  for $0<r\leq 1$, $t\geq 0$ and $w(0,t)=0$ for $t\geq 0$. Similarly as in the
  proof of Lemma 4.1 we can prove the following estimates:
$$
  \max_{0\leq r\leq 1}|w(r,t)-v^*(r)|\leq C\delta r(1-r)e^{-\mu t}
  \quad \mbox{for} \;\; t\geq 0,
\eqno{(4.14)}
$$
$$
  \max_{0\leq r\leq 1}\Big|{\partial w(r,t)\over\partial r}-v^{*'}(r)\Big|
  \leq C\delta e^{-\mu t} \quad \mbox{for} \;\; t\geq 0.
\eqno{(4.15)}
$$
  Indeed, for $0<r\leq 1/2$ we have
$$
\begin{array}{rl}
  |w(r,t)-v^*(r)|\leq
  &\displaystyle\frac{1}{r^2}\int_0^r|g(m(\rho,z(t))+\eta(\rho,t),
   p(\rho,t))-g(m(\rho,z^*),p^*(\rho))|\rho^2d\rho
\\ [0.3cm]
  &\displaystyle+r\int_0^1|g(m(\rho,z(t))+\eta(\rho,t),
   p(\rho,t))-g(m(\rho,z^*),p^*(\rho))|\rho^2d\rho
\\
  \leq &\displaystyle Cr[|z(t)-z^*|
  +\max_{0\leq r\leq 1}|\eta(r,t)|+\max_{0\leq r\leq 1}|p(r,t)-p^*(r)|]
\\
  \leq &Cr\delta e^{-\mu t} \quad \mbox{for}\;\;  t\geq 0,
\end{array}
$$
  and for $1/2<r\leq 1$ we have
$$
\begin{array}{rl}
  |w(r,t)-v^*(r)|\leq
  &\displaystyle\Big(\frac{1}{r^2}-r\Big)\int_0^1|g(m(\rho,z(t))+\eta(\rho,t),
   p(\rho,t))-g(m(\rho,z^*),p^*(\rho))|\rho^2d\rho
\\ [0.3cm]
  &\displaystyle+\frac{1}{r^2}\int_r^1|g(m(\rho,z(t))+\eta(\rho,t),
   p(\rho,t))-g(m(\rho,z^*),p^*(\rho))|\rho^2d\rho
\\
  \leq &\displaystyle C(1-r)[|z(t)-z^*|
  +\max_{0\leq r\leq 1}|\eta(r,t)|+\max_{0\leq r\leq 1}|p(r,t)-p^*(r)|]
\\
  \leq &C(1-r)\delta e^{-\mu t} \quad \mbox{for}\;\;  t\geq 0.
\end{array}
$$
  Combining the above two estimates together we obtain (4.14). The proof of
  (4.15) is simplar; we omit it.

  Having proved (4.14) and (4.15), it follows that all the results in Section 5
  of \cite{Cui3} can be applied to the family of similarity transforms
  $\bar{r}\mapsto r=S(\bar{r};t,s)$ ($t\geq s\geq 0$) of the unit interval
  $[0,1]$ to itself obtained from solving the initial value problem
$$
\left\{
\begin{array}{l}
  \displaystyle {\partial r\over\partial t}+
  v^*(\bar{r}){\partial r\over\partial\bar{r}}=w(r,t)
  \quad \mbox{for} \;\; 0\leq \bar{r}\leq 1,\;\; t>s,\\[0.2cm]
  \displaystyle r|_{t=s}=\bar{r} \quad \mbox{for} \;\; 0\leq \bar{r}\leq 1.
\end{array}
\right.
\eqno{(4.16)}
$$
  The estimates (4.12) and (4.13) then follows from a similar argument as in the
  proof of Lemma 6.4 of \cite{Cui3}. This proves Lemma 4.2. $\quad\Box$
\medskip

  Let $\mu^*$ be the constant as in Lemma 4.2 and fix a number $0<\mu<\mu^*$.
  Later on the notation $\mu$ always denotes this fixed number.
\medskip

  {\bf Lemma 4.3}\ \ {\em  Let $\eta\in C([0,1]\times [0,\infty))$ be given such
  that $\eta_r\in C([0,1]\times [0,\infty))$ and
$$
  \sup_{0\leq r\leq 1}|\eta(r,t)|+\sup_{0\leq r\leq 1}|\eta_r(r,t)|
  \leq\delta e^{-\mu t} \quad \mbox{for} \;\; t\geq 0,
\eqno{(4.17)}
$$
  where $\delta$ and $\mu$ are given positive constants, $\delta$ sufficiently
  small and $\mu<\mu^*$. Let $V_0\in X_0$ be such that $\|V_0\|_{X_0}\leq\delta$.
  Under these conditions, the equation $(4.2)$ subject to the initial condition
  $V(0)=V_0$ has a unique solution $V\in C([0,\infty),X_0)\cap C^1([0,\infty),X)$
  satisfying the following estimates:
$$
  \|V(t)\|_{X_0}\leq C\delta e^{-\mu t}, \quad
  \|\dot{V}(t)\|_X\leq C\delta e^{-\mu t} \quad \mbox{for}\;\; t\geq 0,
\eqno{(4.18)}
$$
  where $C$ is a positive constant independent of $\eta$, $V_0$ and $\delta$.}
\medskip

  {\em Proof}:\ \ Let ${\mathbf M}$ be the set of all functions $V=V(t)\in
  C([0,\infty),X)$ satisfying the following conditions:
$$
  V(0)=V_0,\quad  \|V(t)\|_{X}\leq C\delta e^{-\mu t}\;\;
  \mbox{for}\;\; t\geq 0,
\eqno{(4.19)}
$$
  where $C$ is a positive constant to be specified later. We introduce a metric
  $d$ on ${\mathbf M}$ by defining $d(V_1,V_2)=\sup_{t\geq 0} e^{\mu t}
  \|V_1(t)-V_2(t)\|_X$ for $V_1,V_2\in {\mathbf M}$.
  It is clear that $({\mathbf M},d)$ is a complete metric space. In what
  follows we split into two steps to prove that the equation (4.2) subject
  to the initial condition $V(0)=V_0$ has a unique solution in ${\mathbf M}$
  provided $\delta$ is sufficiently small, and the solution also satisfies the
  other two estimates in (4.19).

  {\em Step 1}:\ \ We prove that if $\delta$ is sufficiently small then for any
  $V\in {\mathbf M}$, the initial value problem
$$
\left\{
\begin{array}{l}
  \dot{U}(t)={\mathbb A}(t,V(t))U(t)+{\mathbb B}(t)U(t)+{\mathbb G}(U(t))+G(t)
  \quad \mbox{for}\;\; t>0,\\
  U(0)=V_0.
\end{array}
\right.
\eqno{(4.20)}
$$
  has a unique solution $U\in C([0,\infty),X_0)\cap C^1([0,\infty),X)$ satisfying
  the following estimates:
$$
  \|U(t)\|_{X}\leq C\delta e^{-\mu t},\quad
  \|U(t)\|_{X_0}\leq C'\delta e^{-\mu t}, \quad
  \|\dot{U}(t)\|_X\leq C'\delta e^{-\mu t} \quad \mbox{for}\;\; t\geq 0,
\eqno{(4.21)}
$$
  where $C$ is the same constant as that appears in (4.19), and $C'$ is another
  positive constant to be specified later. To this end we let
$$
  \widetilde{\mathbf M}=\{U\in C([0,\infty),X_0):
  \|U(t)\|_{X}\leq C\delta e^{-\mu t}\;\;\mbox{and} \;\;
  \|U(t)\|_{X_0}\leq C'\delta e^{-\mu t} \;\; \mbox{for}\;\; t\geq 0\},
$$
  and introduce a metric $d$ on it by defining $d(U_1,U_2)=\sup_{t\geq 0}
  e^{\mu t}\|U_1(t)-U_2(t)\|_{X_0}$ for $U_1,U_2\in \widetilde{\mathbf M}$.
  $(\widetilde{\mathbf M},d)$ is clearly a complete metric space. Given $U\in
  \widetilde{\mathbf M}$, we consider the following initial value problem:
$$
\left\{
\begin{array}{l}
  \displaystyle{d\widetilde{U}(t)\over dt}={\mathbb A}(t,V(t))\widetilde{U}(t)
  +{\mathbb B}(t)U(t)+{\mathbb G}(U(t))+G(t) \quad \mbox{for}\;\; t>0,\\
  \widetilde{U}(0)=V_0.
\end{array}
\right.
$$
  Since $U(t)\in C([0,\infty),X_0)$, we have ${\mathbb G}(U(t))\in
  C([0,\infty),X_0)$. Since also $G\in C([0,\infty),X_0)$, it follows that the
  above problem has a unique solution $\widetilde{U}\in C([0,\infty),X_0)\cap
  C^1([0,\infty),X)$, given by
$$
  \widetilde{U}(t)={\mathbb U}(t,0,V)V_0+
  \int_0^t{\mathbb U}(t,s,V)[{\mathbb B}(s)U(s)+{\mathbb G}(U(s))ds+G(s)]ds.
\eqno{(4.22)}
$$
  Choose a number $\mu'$ such that $\mu<\mu'<\mu^*$. By Lemma 4.2, if
  $\delta$ is sufficiently small then there exist constants $C_1,C_2>0$
  depending only on $\mu'$ such that for any $V=V(t)\in C([0,\infty),X)$
  satisfying (4.19), the following estimates hold:
$$
  \|{\mathbb U}(t,s, V)\|_{L(X)}\leq C_1 e^{-\mu'(t-s)} \quad
  \mbox{for} \;\; t\geq s\geq 0,
\eqno{(4.23)}
$$
$$
  \|{\mathbb U}(t,s, V)\|_{L(X_0)}\leq C_2 e^{-\mu'(t-s)} \quad
  \mbox{for} \;\; t\geq s\geq 0,
\eqno{(4.24)}
$$
  Using (4.22), (4.23), (4.4), (4.8), the fact that $\|V_0\|_X\leq
  \|V_0\|_{X_0}\leq\delta$ and the condition $\|U(t)\|_{X}\leq C\delta
  e^{-\mu t}$ we see that for some positive constant $C''$,
$$
\begin{array}{rcl}
  \|\widetilde{U}(t)\|_{X}&\leq &\displaystyle C_1 e^{-\mu' t}\|V_0\|_{X}+
  C_1\int_0^t e^{-\mu'(t-s)}[\|{\mathbb B}(s)U(s)\|_{X}
  +\|{\mathbb G}(U(s))\|_{X}+\|G(s)\|_{X}]ds
\\ [0.2cm]
  &\leq &\displaystyle C_1\delta e^{-\mu' t}+
  C_1C''\int_0^t e^{-\mu'(t-s)}[C\delta^2e^{-2\mu s}+
  C^2\delta^2e^{-2\mu s}+\delta e^{-\mu s}]ds
\\
  &\leq &\displaystyle\Big(1+2C''C^2\delta\Big)C_1\delta\int_0^t e^{-\mu(t-s)} e^{-2\mu s}ds
  +\frac{C_1C''\delta}{\mu'-\mu}e^{-\mu' t}\int_0^t e^{(\mu'-\mu)s}ds
\\
  &\leq &\displaystyle\Big(1+\frac{2}{\mu}C''C^2\delta
  +\frac{C''}{\mu'-\mu}\Big)C_1\delta e^{-\mu t}.
\end{array}
$$
  Hence, if we first choose the constant $C>0$ sufficiently large such that
  $C\geq\displaystyle\Big(2+\frac{C''}{\mu'-\mu}\Big)C_1$ and next choose
  $\delta_0>0$ sufficiently small such that $\displaystyle\frac{2}{\mu}C''C^2
  \delta_0\leq 1$, then for all $0<\delta\leq\delta_0$ we have
$$
  \|\widetilde{U}(t)\|_{X}\leq C\delta e^{-\mu t}
   \quad \mbox{for all}\;\; t\geq 0.
$$
  Similalrly, by using (4.22), (4.24), (4.5), (4.8), the fact that $\|V_0\|_{X_0}
  \leq\delta$ and the condition $\|U(t)\|_{X_0}\leq C'\delta e^{-\mu t}$ we see
  that by first choosing $C'$ sufficiently large and next choosing $\delta_0$
  further small (when necessary), we also have
$$
  \|\widetilde{U}(t)\|_{X_0}\leq C'\delta e^{-\mu t}
   \quad \mbox{for all}\;\; t\geq 0.
$$
  Hence $\widetilde{U}\in\widetilde{\mathbf M}$. We now define a map
  $\widetilde{\mathbf S}:\widetilde{\mathbf M}\to\widetilde{\mathbf M}$ by setting
  $\widetilde{\mathbf S}(U)=\widetilde{U}$ for every $U\in \widetilde{\mathbf
  M}$. A similar argument by using (4.7) instead of either (4.4) or (4.5)
  we can prove that by choosing $\delta_0$ further small when necessary,
  $\widetilde{\mathbf S}$ is a contraction mapping. It follows from the
  Banach fixed point theorem that $\widetilde{\mathbf S}$ has a fixed point in
  $\widetilde{\mathbf M}$, which is clearly a solution of the problem (4.20)
  in $C([0,\infty),X_0)$. Uniqueness of the solution follows from a standard
  argument. The assertion that the solution $U\in C^1([0,\infty),X_0)$ and
  it satisfies the third estimate in (4.21) are easy consequences of the
  equation in the first line of (4.20) (cf. the proof of Lemma 7.1 of \cite{Cui3}).

  The assertion obtained in Step 1 in particular implies that for every $V$ in
  ${\mathbf M}$, the solution $U$ of (7.2) also belongs to ${\mathbf M}$. Thus
  we can define a mapping ${\mathbf S}:{\mathbf M}\to {\mathbf M}$ as follows:
  For any $V\in {\mathbf M}$,
$$
  {\mathbf S}(V)=U=\mbox{the solution of (4.20)}.
$$

  {\em Step 2}:\ \ We prove that if $\delta$ is sufficiently small, ${\mathbf S}$
  is a contraction mapping. For this purpose, let $V_1,V_2\in {\mathbf M}$ and
  denote $U_1={\mathbf S}(V_1)$, $U_2={\mathbf S}(V_2)$ and $W=U_1-U_2$. Then
  $W$ satisfies
$$
\left\{
\begin{array}{l}
  \displaystyle{dW(t)\over dt}={\mathbb A}(t,V_1(t))W(t)+[{\mathbb A}(t,V_1(t))-
  {\mathbb A}(t,V_2(t))]U_2(t)\\
  \qquad\qquad +{\mathbb B}(t)W(t)+[{\mathbb G}(U_1(t))-{\mathbb G}(U_2(t))]
  \quad \mbox{for}\;\; t>0,\\
  W(0)=0,
\end{array}
\right.
$$
  so that
$$
  W(t)=\int_0^t {\mathbb U}(t,s,V_1)\{[{\mathbb A}(s,V_1(s))-
  {\mathbb A}(s,V_2(s))]U_2(s)ds+{\mathbb B}(s)W(s)+
  [{\mathbb G}(U_1(s))-{\mathbb G}(U_2(s))]\}ds.
$$
  It can be easily shown that (cf. the proof of Lemma 7.2 of \cite{Cui3})
$$
  \|[{\mathbb A}(s,V_1(s))-{\mathbb A}(s,V_2(s))]U_2(s)\|_X
  \leq C\|V_1(s)-V_2(s)\|_X\|U_2(s)\|_{X_0}
  \leq C\delta e^{-2\mu s}d(V_1,V_2).
\eqno{(4.25)}
$$
  Besides, from (4.3) we have
$$
  \|{\mathbb B}(s)W(s)\|_X\leq C\delta e^{-\mu s}\|W(s)\|_X
  \leq C\delta e^{-2\mu s}d(U_1,U_2),
$$
  and from (4.6) we have
$$
\begin{array}{rcl}
  \|{\mathbb G}(U_1(s))-{\mathbb G}(U_2(s))\|_X
  &\leq & C\big(\|U_1(s)\|_X+\|U_2(s)\|_X\big)\|U_1(s)-U_2(s)\|_X
\\
  &\leq & C\delta e^{-\mu s}\|W(s)\|_X\leq C\delta e^{-2\mu s}d(U_1,U_2).
\end{array}
$$
  From these relations and (4.21) we get
$$
\begin{array}{rcl}
  \|U_1(t)-U_2(t)\|_X
  &\leq &\displaystyle C\delta d(V_1,V_2)\int_0^t e^{-\mu (t-s)}e^{-2\mu s}ds
  +C\delta d(U_1,U_2)\int_0^t e^{-\mu (t-s)}e^{-2\mu s}ds
\\
  &\leq &\displaystyle C\delta e^{-\mu t}d(V_1,V_2)+
   C\delta e^{-\mu t}d(U_1,U_2),
\end{array}
$$
  which yields $d(U_1,U_2)\leq C\delta d(V_1,V_2)+ C\delta d(U_1,U_2)$. The
  desired assertion now easily follows.

  It follows that if $\delta$ is sufficiently small then the map ${\mathbf S}$
  has a fixed point $U$ in ${\mathbf M}$. Since the image of ${\mathbf S}$ is
  contained in $\widetilde{\mathbf M}$, we obtain the assertion of Lemma
  4.3. This completes the proof of Lemma 4.3. $\quad\Box$
\medskip

  {\em Remark}:\ \ A similar argument as in Step 2 of the above proof shows
  that if $V_1,V_2$ are solutions of the equation (4.2) with respect to
  $\eta_1,\eta_2\in Y$ respectively (with same initial data), then we have the
  following estimate:
$$
  d(V_1,V_2)\leq C\delta\sup_{0\leq r\leq 1\atop t\geq 0}
  e^{\mu t}|\eta_1(r,t)-\eta_2(r,t)|.
\eqno{(4.26)}
$$
  Indeed, to emphasize the dependence of the operator $\mathbb{A}(t,V)$ on $\eta$
  we redenote it as $\mathbb{A}(\eta,V)$. Similarly as in (4.25) we have
$$
\begin{array}{rcl}
  &&\|[{\mathbb A}(\eta_1(\cdot,s),V_1(s))-{\mathbb A}(\eta_2(\cdot,s),V_2(s))]U_2(s)\|_X
\\
  &\leq &\displaystyle C\sup_{0\leq r\leq 1}|\eta_1(r,s)-\eta_2(r,s)|
  \|U_2(s)\|_{X_0}+ C\|V_1(s)-V_2(s)\|_X\|U_2(s)\|_{X_0}
\\
  &\leq &\displaystyle C e^{-2\mu s}\delta\sup_{0\leq r\leq 1\atop t\geq 0}
  e^{\mu t}|\eta_1(r,t)-\eta_2(r,t)|+C\delta e^{-2\mu s}d(V_1,V_2).
\end{array}
$$
  From this inequality and a similar argument as before we obtain the desired assertion.

  The estimate (4.26) will be useful in Section 6.

\section{Decay estimates of the solution of (2.8)}
\setcounter{equation}{0}

\hskip 2em
  {\bf Lemma 5.1}\ \ {\em Let $(\eta,p,z)=(\eta(r,t),p(r,t),z(t))$ $(0\leq r\leq
  1, t\geq 0)$ be given such that $\eta,\eta_r\in C([0,1]\times[0,\infty))(=
  C([0,\infty),X))$, $p\in C([0,\infty),X_0)$, $z\in C^1[0,\infty)$. Assume that
$$
  \sup_{0\leq r\leq 1}|\eta(r,t)|+\sup_{0\leq r\leq 1}|\eta_r(r,t)|
  \leq\delta e^{-\mu t} \quad \mbox{for} \;\; t\geq 0,
\eqno{(5.1)}
$$
$$
  \sup_{0\leq r\leq 1}|p(r,t)-p^*(r)|+\sup_{0<r<1}|r(1-r)[p_r(r,t)-p^{*'}(r)]|
  \leq\delta e^{-\mu t} \quad \mbox{for} \;\; t\geq 0,
\eqno{(5.2)}
$$
$$
  |z(t)-z^*|+|\dot{z}(t)|\leq\delta e^{-\mu t} \quad \mbox{for} \;\; t\geq 0.
\eqno{(5.3)}
$$
  Let $\tilde{\eta}=\tilde{\eta}(r,t)$ $(0\leq r\leq 1, t\geq 0)$ be the solution
  of (2.8) with initial data $\tilde{\eta}(r,0)=\eta_0(r)$ satisfying
$$
  \sup_{0\leq r\leq 1}|\eta_0(r)|+\sup_{0\leq r\leq 1}|\eta_0'(r)|
  \leq\varepsilon'\delta,
\eqno{(5.4)}
$$
  where $\varepsilon'$ is a given small positive constant. Under these conditions,
  there exist constants $\varepsilon_0,\varepsilon_0'>0$ and $C>0$ independent
  of $\delta$ such that if $0<\varepsilon\leq\varepsilon_0$ and $0<\varepsilon'
  \leq\varepsilon_0'$ then
$$
   |\tilde{\eta}(r,t)|+|\tilde{\eta}_r(r,t)|\leq
    C(\varepsilon+\varepsilon')\delta e^{-\mu t}
\eqno{(5.5)}
$$
   for all $0\leq r\leq 1$ and $t\geq 0$.}
\medskip

  {\em Proof}:\ \ Let $c_0=\displaystyle\min_{-\delta_0\leq c\leq 1+\delta_0}F'(c)
  >0$ and $\lambda=c_0e^{2(z^*-\delta_0)}$. Since $\mathbbm{g}(0,p^*,z^*)=v^*(1)
  =0$, we have
\setcounter{equation}{6}
\begin{eqnarray}
  \mathbbm{g}(\eta(\cdot,t),p(\cdot,t),z(t))
  &=&|\mathbbm{g}(\eta(\cdot,t),p(\cdot,t),z(t))-\mathbbm{g}(0,p^*,z^*)|
\nonumber\\ [0.2cm]
  &\leq &\displaystyle\int_0^1|g(m(\rho,z(t))+\eta(\rho,t),p(\rho,t))
  -g(m(\rho,z^*),p^*(\rho))|\rho^2d\rho
\nonumber\\ [0.1cm]
  &\leq &\displaystyle C[\max_{0\leq r\leq 1}|\eta(r,t)|
  +\max_{0\leq r\leq 1}|p(r,t)-p^*(r)|+|z(t)-z^*|]
\nonumber\\
  &\leq & C\delta e^{-\mu t}.
\end{eqnarray}
  Besides, we also have
$$
\begin{array}{c}
  z^*-\delta_0\leq z^*-\delta e^{-\mu t}\leq z(t)
  \leq z^*+\delta e^{-\mu t}\leq z^*+\delta_0,
\\ [0.1cm]
  |m_z(r;z)-r m_r(r;z)|\leq Ce^{2z(t)}\leq Ce^{2(z^*+\delta_0)},
\\ [0.1cm]
  \displaystyle a(r;\eta(r,t),z(t))=\int_0^1F'(m(r,z(t))+\theta\eta(r,t))d\theta
  \geq c_0.
\end{array}
$$
  It follows that
$$
  |\varepsilon e^{2z(t)}\mathbbm{g}(\eta,p,z)\big[m_z(r;z)-r m_r(r;z)\big]|
  \leq C\varepsilon\delta e^{-\mu t},
$$
  and
$$
  e^{2z(t)}a(r;\eta,z)\geq c_0e^{2(z^*-\delta_0)}=\lambda.
$$
  Using these estimates, we can easily verify that
$$
  \eta_+(r,z)=2\lambda^{-1}C\varepsilon\delta e^{-\mu t}
  +\varepsilon'\delta e^{-\lambda t/\varepsilon}
  \quad \mbox{and} \quad
  \eta_-(r,z)=-2\lambda^{-1}C\varepsilon\delta e^{-\mu t}
  -\varepsilon'\delta e^{-\lambda t/\varepsilon}
$$
  are respectively upper and lower solutions of the problem (2.8) provided
  $\varepsilon$ is sufficiently small. Hence, by the maximum principle we
  see that
$$
  -2\lambda^{-1}C\varepsilon\delta e^{-\mu t}
  -\varepsilon'\delta e^{-\lambda t/\varepsilon}
  \leq\tilde{\eta}(r,z)
  \leq 2\lambda^{-1}C\varepsilon\delta e^{-\mu t}
  +\varepsilon'\delta e^{-\lambda t/\varepsilon}
\eqno{(5.8)}
$$
  for all $0\leq r\leq 1$ and $t\geq 0$.

  Next, by differentiating the equation in the first line of (2.8) in $r$
  and using the boundary value conditions for $\tilde{\eta}$, we see that
  $\tilde{\eta}_r$ satisfies
$$
\left\{
\begin{array}{l}
   \displaystyle\varepsilon e^{2z(t)}(\tilde{\eta}_r)_t=(\tilde{\eta}_r)_{rr}
  +\big[{2\over r}+\varepsilon e^{2z(t)}r\mathbbm{g}(\eta,p,z)\big](\tilde{\eta}_r)_r
  -\big[e^{2z(t)}a(r;\eta,z)+{2\over r^2}-\varepsilon e^{2z(t)}\mathbbm{g}(\eta,p,z)\big]
  \tilde{\eta}_r
\\ [0.2cm]
   \quad\quad\quad\quad -e^{2z(t)}a_r(r;\eta,z)\tilde{\eta}
   -e^{2z(t)}a_\eta(r;\eta,z)\tilde{\eta}\eta_r
\\
   \quad\quad\quad\quad -\varepsilon e^{2z(t)}\mathbbm{g}(\eta,p,z)\big[m_{rz}(r;z)
   -rm_{rr}(r;z)-m_r(r;z)\big]
   \quad \mbox{for}\;\; 0<r<1,\;\; t>0,
\\
   \tilde{\eta}_r(0,t)=0, \quad \tilde{\eta}_{rr}(1,t)
   +[2+\varepsilon e^{2z(t)}\mathbbm{g}(\eta,p,z)]\tilde{\eta}_r(1,t)=
   -\varepsilon e^{2z(t)}\mathbbm{g}(\eta,p,z)m_r(1;z)\;\;  \mbox{for}\;\; t>0,
\\
   \tilde{\eta}_r(r,0)=\eta_0'(r) \quad \mbox{for} \;\; 0\leq r\leq 1.
\end{array}
\right.
\eqno{(5.9)}
$$
  If $\varepsilon$ is sufficiently small we have
$$
  e^{2z(t)}a(r;\eta,z)+{2\over r}-\varepsilon e^{2z(t)}\mathbbm{g}(\eta,p,z)\geq
  e^{2(z^*-\delta_0)}c_0+2-C\varepsilon\delta_0>\lambda.
$$
  Moreover, using the assumptions on $\eta,z$ and also using (5.7) and (5.8),
  we see that
$$
\begin{array}{rl}
  &|e^{2z(t)}a_r(r;\eta,z)\tilde{\eta}
   +e^{2z(t)}a_\eta(r;\eta,z)\tilde{\eta}\eta_r
   +\varepsilon e^{2z(t)}\mathbbm{g}(\eta,p,z)\big[m_{rz}(r;z)
   -rm_{rr}(r;z)-m_r(r;z)\big]|
\\ [0.3cm]
  \leq & C\varepsilon\delta e^{-\mu t}
  +C\varepsilon'\delta e^{-\lambda t/\varepsilon}.
\end{array}
$$
  Since for $\varepsilon$ sufficiently small we also have
$$
  2+\varepsilon e^{2z(t)}\mathbbm{g}(\eta,p,z)\geq 1
  \quad \mbox{and} \quad
  |\varepsilon e^{2z(t)}\mathbbm{g}(\eta,p,z)\big[m_z(1;z)-m_r(1;z)\big]|
  \leq C\delta\varepsilon e^{-\mu t},
$$
  again by using the maximum principle we conclude that if $\varepsilon$ and
  $\varepsilon'$ are sufficiently small then
$$
  -(1+2\lambda^{-1})C\varepsilon\delta e^{-\mu t}
  -\varepsilon'\delta e^{-\lambda t/\varepsilon}
  \leq\tilde{\eta}_r(r,z)
  \leq (1+2\lambda^{-1})C\varepsilon\delta e^{-\mu t}
  +C\varepsilon'\delta e^{-\lambda t/\varepsilon}
\eqno{(5.10)}
$$
  for all $0\leq r\leq 1$ and $t\geq 0$. Combining (5.8) and (5.10), and
  choosing $\varepsilon_0$ further small when necessary so that
  $\lambda/\varepsilon\geq\mu$ for $0<\varepsilon\leq\varepsilon_0$, we
  obtain (5.5). This proves Lemma 5.1. $\quad\Box$
\medskip

  {\bf Lemma 5.2}\ \ {\em Under the same conditions as in Lemma 5.1, there also
  holds the following estimate:
$$
   \int_0^t\!\!\int_0^1e^{-\frac{2c_0}{\varepsilon}(t-s)}
   |\tilde{\eta}_{rr}(r,s)|^2r^2drds\leq C\varepsilon\delta^2 e^{-2\mu t}
\eqno{(5.11)}
$$
   for all $t\geq 0$.}
\medskip

  {\em Proof}:\ \ From the equation in the first three lines in (5.9) we have
$$
\begin{array}{l}
   \displaystyle\varepsilon(\tilde{\eta}_r)_t= e^{-2z(t)}r^{-2}(r^2\tilde{\eta}_{rr})_r
  +\varepsilon r\mathbbm{g}(\eta,p,z)\tilde{\eta}_{rr}
  -\big[a(r;\eta,z)+{2\over r^2}e^{-2z(t)}-\varepsilon\mathbbm{g}(\eta,p,z)\big]
  \tilde{\eta}_r
\\ [0.2cm]
   \quad\quad\quad\quad -a_r(r;\eta,z)\tilde{\eta}-a_\eta(r;\eta,z)\tilde{\eta}\eta_r
  -\varepsilon\mathbbm{g}(\eta,p,z)\big[m_{rz}(r;z)
   -rm_{rr}(r;z)-m_r(r;z)\big].
\end{array}
$$
  Multiply both sides of this equation with $\tilde{\eta}_rr^2$ and next
  integrating with respect to $r$ in $[0,1]$, we get
$$
\begin{array}{rl}
  &\displaystyle\frac{\varepsilon}{2}\frac{d}{dt}
  \int_0^1|\tilde{\eta}_r|^2r^2dr
\\ [0.2cm]
  =&\displaystyle-e^{-2z(t)}\int_0^1|\tilde{\eta}_{rr}|^2r^2dr
  +\varepsilon\mathbbm{g}(\eta,p,z)\int_0^1\tilde{\eta}_{rr}\tilde{\eta}_rr^3dr
  -\int_0^1a(r;\eta,z)|\tilde{\eta}_r|^2r^2dr
\\ [0.2cm]
  &\displaystyle-2e^{-2z(t)}\int_0^1|\tilde{\eta}_r|^2dr
  +\varepsilon\mathbbm{g}(\eta,p,z)\int_0^1|\tilde{\eta}_r|^2r^2dr
  -\int_0^1a_r(r;\eta,z)\tilde{\eta}\tilde{\eta}_rr^2dr
\\ [0.2cm]
  &\displaystyle-\int_0^1a_\eta(r;\eta,z)\tilde{\eta}\eta_r\tilde{\eta}_rr^2dr
  -\varepsilon\mathbbm{g}(\eta,p,z)\int_0^1\big[m_{rz}(r;z)
   -rm_{rr}(r;z)-m_r(r;z)\big]\tilde{\eta}_rr^2dr.
\end{array}
$$
  We have
$$
  e^{-2z(t)}\int_0^1|\tilde{\eta}_{rr}|^2r^2dr\geq
  e^{-2(z^*+\delta_0)}\int_0^1|\tilde{\eta}_{rr}|^2r^2dr
  \stackrel{{\rm def}}{=}C_0\int_0^1|\tilde{\eta}_{rr}|^2r^2dr,
$$
$$
\begin{array}{rl}
   \displaystyle\varepsilon\mathbbm{g}(\eta,p,z)
   \int_0^1\tilde{\eta}_{rr}\tilde{\eta}_rr^3dr
   =&\displaystyle\frac{\varepsilon}{2}
   \mathbbm{g}(\eta,p,z)|\tilde{\eta}_r(1,t)|^2
  -\frac{3\varepsilon}{2}\mathbbm{g}(\eta,p,z)
  \int_0^1|\tilde{\eta}_r|^2r^2dr
\\ [0.2cm]
   \leq&\displaystyle C\varepsilon^3\delta^3e^{-3\mu t}
   \quad (\mbox{by (5.5) and (5.7)}),
\end{array}
$$
$$
  \int_0^1a(r;\eta,z)|\tilde{\eta}_r|^2r^2dr\geq
  c_0\int_0^1|\tilde{\eta}_r|^2r^2dr,
$$
$$
  2e^{-2z(t)}\int_0^1|\tilde{\eta}_r|^2dr\geq 0.
$$
  Moreover, by using (5.5), (5.7) and the boundedness of $\eta,p,z$ (with
  bounds independent of special choice of these functions) we have
$$
\begin{array}{rl}
   &\displaystyle\Big|\varepsilon\mathbbm{g}(\eta,p,z)
   \int_0^1|\tilde{\eta}_r|^2r^2dr
  -\int_0^1a_r(r;\eta,z)\tilde{\eta}\tilde{\eta}_rr^2dr
  -\int_0^1a_\eta(r;\eta,z)\tilde{\eta}\eta_r\tilde{\eta}_rr^2dr
\\ [0.2cm]
  &\displaystyle-\varepsilon\mathbbm{g}(\eta,p,z)\int_0^1\big[m_{rz}(r;z)
   -rm_{rr}(r;z)-m_r(r;z)\big]\tilde{\eta}_rr^2dr\Big|
\\ [0.2cm]
   \leq &\displaystyle C\varepsilon^3\delta^3e^{-3\mu t}
   +C\varepsilon^2\delta^2e^{-2\mu t}
   +C\varepsilon^2\delta^3e^{-3\mu t}
   +C\varepsilon^2\delta^2e^{-2\mu t}
   \leq C\varepsilon^2\delta^2e^{-2\mu t}.
\end{array}
$$
  It follows that
$$
  \frac{\varepsilon}{2}\frac{d}{dt}\int_0^1|\tilde{\eta}_r|^2r^2dr
  \leq -C_0\int_0^1|\tilde{\eta}_{rr}|^2r^2dr
  -c_0\int_0^1|\tilde{\eta}_r|^2r^2dr+C\varepsilon^2\delta^2e^{-2\mu t}.
$$
  From this inequality we easily deduce that if $\varepsilon$ is sufficiently
  small such that $0<\varepsilon<\frac{c_0}{2\mu}$, then
$$
\begin{array}{rl}
   &\displaystyle\int_0^1|\tilde{\eta}_r|^2r^2dr+\frac{2C_0}{\varepsilon}
  \int_0^t\!\!\int_0^1e^{-\frac{2c_0}{\varepsilon}(t-s)}
   |\tilde{\eta}_{rr}(r,s)|^2r^2drds
\\ [0.2cm]
   \leq &\displaystyle e^{-\frac{2c_0}{\varepsilon}t}\int_0^1|\eta_0'(r)|^2r^2dr
   +C\varepsilon^2\delta^2e^{-2\mu t}
   \leq\displaystyle C\delta^2e^{-4\mu t}+C\varepsilon^2\delta^2e^{-2\mu t}.
\end{array}
$$
  The estimate (5.11) now immediately follows. $\quad\Box$
\medskip

  {\bf Corollary 5.3}\ \ {\em Under the same conditions as in Lemma 5.1, for any
  $0<\nu<\mu$ there exists corresponding constant $C>0$ such that the following
  estimate holds:
$$
   \int_0^{\infty}\!\!\int_0^1e^{2\nu t}|\tilde{\eta}_{rr}(r,t)|^2r^2drdt
   \leq C\delta^2.
\eqno{(5.12)}
$$
}

  {\em Proof}:\ \ Multiplying both sides of (5.11) with $e^{2\nu t}$ and integrating
  with respect to $t$ over $[0,\infty)$, we get
$$
   \int_0^{\infty}\!\!\int_0^t\!\!\int_0^1e^{2\nu t-\frac{2c_0}{\varepsilon}(t-s)}
   |\tilde{\eta}_{rr}(r,s)|^2r^2drdsdt
   \leq C\varepsilon\delta^2\int_0^{\infty}\! e^{-2(\mu-\nu)t}dt
   =\frac{C\varepsilon\delta^2}{2(\mu-\nu)}.
$$
  Since
$$
\begin{array}{rl}
   &\displaystyle\int_0^{\infty}\!\!\int_0^t\!\!\int_0^1
   e^{2\nu t-\frac{2c_0}{\varepsilon}(t-s)}
   |\tilde{\eta}_{rr}(r,s)|^2r^2drdsdt
\\ [0.3cm]
   = &\displaystyle\int_0^{\infty}\!\!\int_0^1\!
   \Big(\!\int_s^{\infty}\!e^{-2(\frac{c_0}{\varepsilon}-\nu)t}dt\Big)
   |\tilde{\eta}_{rr}(r,s)|^2r^2drds
\\ [0.3cm]
   = &\displaystyle\frac{\varepsilon}{2(c_0-\nu\varepsilon)}
   \int_0^{\infty}\!\!\int_0^1\!|\tilde{\eta}_{rr}(r,s)|^2r^2drds,
\end{array}
$$
  we see that (5.12) follows. $\quad\Box$
\medskip

  {\bf Corollary 5.4}\ \ {\em Under the same conditions as in Lemma 5.1, for any
  $0<\nu<\mu$ there exists corresponding constant $C>0$ such that the following
  estimate holds:
$$
   \int_0^{\infty}\!\!\int_0^1e^{2\nu t}|\tilde{\eta}_t(r,t)|^2r^2drdt
   \leq\frac{C\delta^2}{\varepsilon}.
\eqno{(5.13)}
$$
}

  {\em Proof}:\ \ This follows from Corollary 5.3 combined with the equation (2.7)
  and the fact that $\displaystyle\max_{0\leq r\leq 1}|\tilde{\eta}(r,t)|\leq
  C\delta e^{-\mu t}$ and $\displaystyle\max_{0\leq r\leq 1}|\tilde{\eta}_r(r,t)|
  \leq C\delta e^{-\mu t}$ for all $t\geq 0$. $\quad\Box$

\section{The proof of Theorem 1.1}

\hskip 2em
  {\em The proof of Theorem 1.1}:\ \ Let $\mu^*$ be the positive constant
  specified in Lemma 4.2, and arbitrarily choose a positive constant $\mu$ such
  that $0<\mu<\mu^*$ and fix it. Let $\delta$ and $\varepsilon'$ be positive
  constants which we shall specify later. We assume that the initial
  data $(c_0(r),p_0(r),z_0)$ of $(c(r,t),p(r,t),z(t))$ satisfy (1.21), (1.22)
  and the following conditions:
\setcounter{equation}{0}
\begin{equation}
  \max_{0\leq r\leq 1}|c_0(r)-c^*(r)|<\varepsilon'\delta,
  \quad
  \sup_{0\leq r\leq 1}|c_0'(r)-c^{*'}(r)|<\varepsilon'\delta,
\end{equation}
\begin{equation}
  \max_{0\leq r\leq 1}|p_0(r)-p^*(r)|<\delta,
  \quad
  \sup_{0<r\leq 1}r(1\!-\!r)|p_0'(r)-p^{*'}(r)|<\delta,
  \quad \mbox{and} \quad
  |z_0-z^*|<\varepsilon'\delta.
\end{equation}
  Let $(Y,\|\cdot\|_Y)$ be the Banach space introduced in Section 2, and consider
  the following set $S\subseteq Y$:
$$
  S=\{\eta\in Y:\;\eta(r,0)=\eta_0(r),\;\|\eta\|_Y\leq\delta\},
$$
  where $\eta_0(r)=c_0(r)-m(r,z_0)$. Let $\eta\in S$ and $V_0=(p_0,z_0)$. It is
  clear that $V_0\in X_0$ and, due to the conditions in (6.2), we have
  $\|V_0\|_{X_0}\leq\delta$. It follows by Lemma 4.3 that there exists positive
  constant $\delta_0$ depending only on the choice of $\mu$ such that if
  $0<\delta\leq\delta_0$ then the equation (4.2) subject to the initial value
  condition $V(0)=V_0$ has a unique solution $V\in C([0,\infty),X_0)\cap
  C^1([0,\infty),X)$ satisfying the estimates in (4.18). Let $V(t)=
  (\varphi(\cdot,t),\zeta(t))$ and set $p(r,t)=p^*(r)+\varphi(r,t)$, $z(t)=z^*+
  \zeta(t)$. From the deduction in Section 4 we see that $(p,z)$ is a
  solution of the system (2.7). Clearly, $(p,z)$ satisfies (1.26) and (1.27), so
  that it also satisfies (5.2) and (5.3) when we replace $\delta$ in those
  conditions with $C\delta$. Now we consider the problem (2.8). It is clear that
  the condition (5.1) is satisfied and, by (6.1) and (6.2), the condition (5.4)
  is also satisfied. It follows by Lemma 5.1 that there exist constants
  $\varepsilon_0,\varepsilon_0'>0$ independent of $\delta$ such that if $0<
  \varepsilon\leq\varepsilon_0$ and $0<\varepsilon'\leq\varepsilon_0'$ then the
  solution $\tilde{\eta}=\tilde{\eta}(r,t)$ of (2.8) satisfies the estimate
  (5.5), so that
$$
   |\tilde{\eta}(r,t)|+|\tilde{\eta}_r(r,t)|
   \leq C(\varepsilon+\varepsilon')\delta e^{-\mu t}
   \leq \delta e^{-\mu t}
$$
  provided $C(\varepsilon_0+\varepsilon_0')\leq 1$. This implies that
  $\tilde{\eta}\in S$ if $\varepsilon$ and $\varepsilon'$ are sufficiently small.
  Let $\mathscr{F}:S\rightarrow S$ be the map $\mathscr{F}(\eta)=\tilde{\eta}$.
  In what follows we prove that $\mathscr{F}$ has a fixed point.

  We first prove that $\mathscr{F}$ is continuous with respect to the metric
$$
   d_1(\eta_1,\eta_2)=\sup_{t\geq 0}e^{\mu t}\Big(
  \int_0^1|\eta_1(r,t)-\eta_2(r,t)|^2r^2dr\Big)^{\frac{1}{2}}.
$$
  Given $\eta_1,\eta_2\in S$, we let $\tilde{\eta}_i=\mathscr{F}(\eta_i)$, $i=1,2$.
  Then $\tilde{\eta}_i$ is the solution of the following problem:
$$
\left\{
\begin{array}{l}
   \displaystyle\varepsilon\tilde{\eta}_{it}=
  e^{-2z_i(t)}r^{-2}(r^2\tilde{\eta}_{ir})_r
  +\varepsilon r\mathbbm{g}(\eta_i,p_i,z_i)\tilde{\eta}_{ir}
  -a(r;\eta_i,z_i)\tilde{\eta}_i
\\
   \quad\quad\quad -\varepsilon\mathbbm{g}(\eta_i,p_i,z_i)
   \big[m_z(r;z_i)-r m_r(r;z_i)\big]
   \quad \mbox{for}\;\; 0<r<1,\;\; t>0,
\\
   \displaystyle\tilde{\eta}_{ir}(0,t)=0,
   \quad \tilde{\eta}_i(1,t)=0 \quad \mbox{for}\;\; t>0,
\\
   \tilde{\eta}_i(r,0)=\eta_0(r) \quad \mbox{for} \;\; 0\leq r\leq 1,
\end{array}
\right.
$$
  where $(p_i,z_i)$ is the solution of the problem (2.7) for $\eta=\eta_i$,
  $i=1,2$. It follows that
$$
\begin{array}{rl}
   &\displaystyle\frac{\varepsilon}{2}\frac{d}{dt}
   \int_0^1|\tilde{\eta}_1-\tilde{\eta}_2|^2r^2dr
\\ [0.2cm]
  =&\displaystyle-e^{-2z_1(t)}\int_0^1|\tilde{\eta}_{1r}-\tilde{\eta}_{2r}|^2r^2dr
  -[e^{-2z_1(t)}-e^{-2z_2(t)}]\int_0^1\tilde{\eta}_{2r}
  (\tilde{\eta}_{1r}-\tilde{\eta}_{2r})r^2dr
\\ [0.2cm]
  &\displaystyle-\frac{3\varepsilon}{2}\mathbbm{g}(\eta_1,p_1,z_1)
   \int_0^1|\tilde{\eta}_1-\tilde{\eta}_2|^2r^2dr
  +\varepsilon [\mathbbm{g}(\eta_1,p_1,z_1)-\mathbbm{g}(\eta_2,p_2,z_2)]
  \int_0^1[\tilde{\eta}_1-\tilde{\eta}_2]\tilde{\eta}_{2r}r^3dr
\\ [0.2cm]
  &\displaystyle-\int_0^1a(r;\eta_1,z_1)|\tilde{\eta}_1-\tilde{\eta}_2|^2r^2dr
  -\int_0^1[a(r;\eta_1,z_1)-a(r;\eta_2,z_2)][\tilde{\eta}_1-\tilde{\eta}_2]
  \tilde{\eta}_2r^2dr
\\ [0.2cm]
   &\displaystyle-\varepsilon\int_0^1\{\mathbbm{g}(\eta_1,p_1,z_1)
   \big[m_z(r;z_1)-r m_r(r;z_1)\big]-\mathbbm{g}(\eta_2,p_2,z_2)
   \big[m_z(r;z_2)-r m_r(r;z_2)\big]\}
\\ [0.2cm]
   &\quad\displaystyle \times[\tilde{\eta}_1-\tilde{\eta}_2]r^2dr.
\end{array}
$$
  Since $a(r;\eta_1,z_1)\geq c_0>0$ and $\eta_i$, $p_i$, $z_i$ ($i=1,2$)
  are bounded functions with bounds independent of special choice of these
  functions, by using (5.5), (5.7) and some standard arguments we see that
  if $\delta$ is sufficiently small then
$$
\begin{array}{rl}
   &\displaystyle\frac{\varepsilon}{2}\frac{d}{dt}
   \int_0^1|\tilde{\eta}_1-\tilde{\eta}_2|^2r^2dr
\\ [0.2cm]
  \leq&\displaystyle-\frac{1}{2}e^{-2z_1(t)}
  \int_0^1|\tilde{\eta}_{1r}-\tilde{\eta}_{2r}|^2r^2dr
  -\frac{1}{2}c_0\int_0^1|\tilde{\eta}_1-\tilde{\eta}_2|^2r^2dr
\\ [0.2cm]
  &\displaystyle+C\delta^2e^{-2\mu t}
  [\max_{0\leq r\leq 1}|\eta_1(r,t)-\eta_2(r,t)|^2
  +\max_{0\leq r\leq 1}|p_1(r,t)-p_2(r,t)|^2+|z_1(t)-z_2(t)|^2]
\\ [0.2cm]
  &\displaystyle+C\varepsilon^2[\max_{0\leq r\leq 1}|\eta_1(r,t)-\eta_2(r,t)|^2
  +\max_{0\leq r\leq 1}|p_1(r,t)-p_2(r,t)|^2+|z_1(t)-z_2(t)|^2].
\end{array}
$$
  Neglecting the first term on the right-hand side of the above inequality,
  assuming that $\varepsilon$ is sufficiently small so that $C\varepsilon\leq
  \frac{1}{2}c_0$, and using the estimate (4.26), we get
$$
  \frac{\varepsilon}{2}\frac{d}{dt}
   \int_0^1|\tilde{\eta}_1-\tilde{\eta}_2|^2r^2dr
   \leq -\frac{c_0}{2}\int_0^1|\tilde{\eta}_1-\tilde{\eta}_2|^2r^2dr
   +C(\delta^2+\varepsilon^2)e^{-2\mu t}\sup_{0\leq r\leq 1\atop t\geq 0}
   e^{2\mu t}|\eta_1(r,t)-\eta_2(r,t)|^2.
$$
  It follows that by assuming that $\varepsilon$ is further small such that
  $c_0-2\mu\varepsilon\geq\frac{1}{2}c_0$ when necessary, we have
$$
   \int_0^1|\tilde{\eta}_1(r,t)-\tilde{\eta}_2(r,t)|^2r^2dr
   \leq C(\delta^2+\varepsilon^2)e^{-2\mu t}
   \Big(\sup_{0\leq r\leq 1\atop t\geq 0}e^{\mu t}|\eta_1(r,t)-\eta_2(r,t)|\Big)^2.
\eqno{(6.3)}
$$
  By applying the three-dimensional interpolation inequality
$$
  \|u\|_{\infty}\leq C\|\nabla u\|_{\infty}^{\frac{3}{5}}
  \|u\|_2^{\frac{2}{5}} \quad \mbox{for}\;\; u\in H_0^1(B_1(0)), \;\;
  \nabla u\in L^{\infty}(B_1(0))
$$
  to the case $u(x)=u(|x|)$, we get the following inequality:
$$
  \sup_{0\leq r\leq 1}|u(r)|\leq C\Big(\sup_{0\leq r\leq 1}|u'(r)|\Big)^{\frac{3}{5}}
  \Big(\int_0^1|u(r)|^2r^2dr\Big)^{\frac{1}{5}}.
\eqno{(6.4)}
$$
  It follows that for any $\eta_1,\eta_2\in S$ we have
$$
  \sup_{0\leq r\leq 1\atop t\geq 0}e^{\mu t}|\eta_1(r,t)-\eta_2(r,t)|\leq
  C\delta^{\frac{3}{5}}\sup_{t\geq 0}\Big(e^{2\mu t}
  \int_0^1|\eta_1(r,t)-\eta_2(r,t)|^2r^2dr\Big)^{\frac{1}{5}}.
\eqno{(6.5)}
$$
  Substituting (6.5) into (6.3) we get
$$
   d_1(\tilde{\eta}_1,\tilde{\eta})\leq C\delta^{\frac{3}{5}}
   (\delta^2+\varepsilon^2)^{\frac{1}{2}}d_1^{\frac{2}{5}}(\eta_1,\eta_2).
$$
  Hence, the map $\mathscr{F}:S\to S$ is continuous with respect to the metric $d_1$.

  Let
$$
\begin{array}{rl}
  S_0=&\Big\{\eta\in S:\;\eta(r,t)\;\mbox{is twice weakly differentiable in}\; r\;
  \mbox{and weakly differentiable in}\; t,\;
\\
  &\;\;\displaystyle\int_0^{\infty}\!\!\int_0^1e^{\mu t}|\eta_{rr}(r,t)|^2r^2drdt
   \leq C\delta^2 \;\; \mbox{and} \; \int_0^{\infty}\!\!\int_0^1e^{\mu t}|\eta_t(r,t)|^2r^2drdt
   \leq\frac{C\delta^2}{\varepsilon}\Big\},
\end{array}
$$
  where $C$ is the constant appearing in (5.12) and (5.13) for the case $\nu=
  \mu/2$. By Lemma 5.3 and Lemma 5.4 we see that $\mathscr{F}(S)\subseteq S_0$.

  Let $\bar{S}$ be the closure of $S$ with respect to the metric $d_1$, and let
  $\bar{\mathscr{F}}:\bar{S}\to\bar{S}$ be the unique continuous extension of the
  map $\mathscr{F}$ onto $\bar{S}$. Then $\bar{\mathscr{F}}(\bar{S})\subseteq
  \bar{S}_0$, where $\bar{S}_0$ is the closure of $S_0$ with respect to the metric
  $d_1$. By using a standard $*$-weak compactness argument we easily see that
$$
\begin{array}{rl}
  \bar{S}_0\subseteq&\Big\{\eta\in C([0,\infty),L^2((0,1), r^2dr))\cap
  L^2((0,\infty),H^2((0,1), r^2dr),e^{\mu t}dt)
\\
  &\;\;\;\;\cap H^1((0,\infty),L^2((0,1), r^2dr), e^{\mu t}dt):\;\eta(r,0)=\eta_0(r)\;
  \; \mbox{for a. e.}\; r\in (0,1),
\\
  &\;\;\;\; \eta(1,t)=0\; \; \mbox{for a. e.}\; t>0,\;\displaystyle{\rm ess}\!\!\!
  \sup_{0\leq r\leq 1\atop t\geq 0}e^{\mu t}|\eta(r,t)|\leq\delta,\;
  \;{\rm ess}\!\!\!\sup_{0\leq r\leq 1\atop t\geq 0}e^{\mu t}|\eta_r(r,t)|\leq\delta,
\\
  &\;\;\displaystyle\;\;
  \int_0^{\infty}\!\!\int_0^1e^{\mu t}|\eta_{rr}(r,t)|^2r^2drdt\leq C\delta^2
  \;\; \mbox{and} \; \int_0^{\infty}\!\!\int_0^1e^{\mu t}|\eta_t(r,t)|^2r^2drdt
   \leq\frac{C\delta^2}{\varepsilon}\Big\}.
\end{array}
$$
  Let
$$
\begin{array}{c}
  \bar{Y}=\displaystyle\Big\{\eta\in C([0,\infty),L^2((0,1), r^2dr)):\,
  \|\eta\|_{\bar{Y}}\stackrel{\mathrm{def}}{=}
  \sup_{t\geq 0}e^{\mu t}\Big(\int_0^1|\eta(r,t)|^2r^2dr\Big)^{\frac{1}{2}}
  <\infty\Big\}.
\end{array}
$$
  It is clear that $(\bar{Y},\|\cdot\|_{\bar{Y}})$ is a Banach space and $\bar{S}$
  is a closed convex subset of this space. For any $T>0$, we let $Y_T$, $\bar{Y}_T$,
  $S_T$, $\bar{S}_T$, $S_{0T}$ and $\bar{S}_{0T}$ be respectively the restrictions
  of $Y$, $\bar{Y}$, $S$, $\bar{S}$, $S_0$ and $\bar{S}_0$ on $[0,1]\times [0,T]$,
  and correspondingly define $\mathscr{F}_T$,
  $\bar{\mathscr{F}}_T$ to be respectively the ``restrictions'' of $\mathscr{F}$,
  $\bar{\mathscr{F}}$ on $S_T$ and $\bar{S}_T$, i.e., for any $\zeta\in S_T$ we
  choose an $\eta\in S$ such that $\eta|_{[0,T]}=\zeta$ and define $\mathscr{F}_T
  (\zeta)=\mathscr{F}(\eta)|_{[0,T]}$, and similarly for $\bar{\mathscr{F}}_T$. By
  uniqueness of solutions of the problems (2.7) and (2.8) in any interval $[0,T]$
  we see that these definitions make sense. Moreover, it is clear that
  $\mathscr{F}_T$ is a continuous self-mapping in $S_T$, $\bar{\mathscr{F}}_T$
  is a continuous self-mapping in $\bar{S}_T$, $\mathscr{F}_T(S_T)\subseteq S_{0T}$
  and $\bar{\mathscr{F}}_T(\bar{S}_T)\subseteq\bar{S}_{0T}$. Note that, in
  particular,
$$
\begin{array}{c}
  \bar{Y}_T=\displaystyle\Big\{\eta\in C([0,T],L^2((0,1), r^2dr)):\,
  \|\eta\|_{\bar{Y}_T}\stackrel{\mathrm{def}}{=}
  \sup_{0\leq t\leq T}e^{\mu t}\Big(\int_0^1|\eta(r,t)|^2r^2dr\Big)^{\frac{1}{2}}
  <\infty\Big\},
\end{array}
$$
$$
\begin{array}{rl}
  \bar{S}_{0T}&\subseteq\Big\{\eta\in C([0,T],L^2((0,1), r^2dr))\cap
  L^2((0,T),H^2((0,1), r^2dr))\cap H^1((0,T),L^2((0,1), r^2dr)):
\\
  &\eta(r,0)=\eta_0(r)\;\; \mbox{for a. e.}\; r\in (0,1),\;
  \eta(1,t)=0\; \; \mbox{for a. e.}\; t\in (0,T),\;\displaystyle
  \sup_{0\leq r\leq 1\atop 0<t<T}e^{\mu t}|\eta(r,t)|\leq\delta,
\\
  &\displaystyle\sup_{0\leq r\leq 1\atop 0<t<T}
  e^{\mu t}|\eta_r(r,t)|\leq\delta,\;\;
  \int_0^T\!\!\int_0^1e^{\mu t}|\eta_{rr}(r,t)|^2r^2drdt\leq C\delta^2
  \;\; \mbox{and}\;\displaystyle\int_0^T\!\!\int_0^1e^{\mu t}|\eta_t(r,t)|^2r^2drdt
  \leq\frac{C\delta^2}{\varepsilon}\Big\}.
\end{array}
$$
  Clearly, $\bar{S}_T$ is a closed convex subset of $\bar{Y}_T$. Moreover,
  it is also easy to see that $\bar{S}_{0T}$ is a compact subset of $\bar{Y}_T$.
  Indeed, let $\{\eta_n\}_{n=1}^{\infty}$ be a sequence in $\bar{S}_{0T}$.
  By compactness of the embedding $H^1((0,1)\times(0,T),r^2drdt)\hookrightarrow
  L^2([0,1]\times[0,T],r^2drdt)$, $*$-weak compactness of bounded sets in
  $L^{\infty}([0,1]\times[0,T])$, and weak compactness of bounded sets in
  $L^2([0,1]\times[0,T])$, we see that there exists a subsequence of
  $\{\eta_n\}_{n=1}^{\infty}$ which we assume, for simplicity of the
  notation, to be $\{\eta_n\}_{n=1}^{\infty}$ itself, and $\eta\in\bar{S}_{0T}$
  such that as $n\to\infty$,
$$
  \eta_n\to\eta \quad \mbox{strongly in}\;\; L^2([0,1]\times[0,T],r^2drdt),
$$
$$
  \eta_{rn}\to\eta_r \quad \mbox{$*$-weakly in}\;\; L^{\infty}([0,1]\times[0,T]),
$$
$$
  \eta_{rrn}\to\eta_{rr} \quad \mbox{weakly in}\;\; L^2([0,1]\times[0,T],r^2drdt),
$$
$$
  \eta_{tn}\to\eta_t \quad \mbox{weakly in}\;\; L^2([0,1]\times[0,T],r^2drdt).
$$
  Integrating the equation
$$
  \frac{d}{dt}\int_0^1|\eta_n(r,t)-\eta(r,t)|^2r^2dr
  =2\int_0^1[\eta_n(r,t)-\eta(r,t)][\eta_{nt}(r,t)-\eta_t(r,t)]r^2dr,
$$
  and using the Cauchy-Schwartz inequality and the inequalities $\displaystyle
  \int_0^T\!\!\int_0^1e^{\mu t}|\eta_{nt}(r,t)|^2r^2drdt
  \leq\frac{C\delta^2}{\varepsilon}$ ($n=1,2,\cdots$) and $\displaystyle
  \int_0^T\!\!\int_0^1e^{\mu t}|\eta_t(r,t)|^2r^2drdt
  \leq\frac{C\delta^2}{\varepsilon}$, we get
$$
  \sup_{0\leq t\leq T}\int_0^1|\eta_n(r,t)-\eta(r,t)|^2r^2dr
  =\frac{C\delta}{\sqrt{\varepsilon}}
  \Big(\int_0^T\!\!\int_0^1|\eta_n(r,t)-\eta(r,t)|^2]r^2drdt\Big)^{\frac{1}{2}},
  \quad n=1,2,\cdots.
$$
  It follows that $\displaystyle\lim_{n\to\infty}\|\eta_n-\eta\|_{\bar{Y}_T}=0$.
  This proves the desired assertion.

  It follows by the Schauder fixed point theorem that $\bar{\mathscr{F}}_T$
  has a fixed point in $\bar{S}_T$, which we denote as $\eta^T$.

  We now prove that $\eta^T\in S_{0T}$ and it is a fixed point of $\mathscr{F}_T$.
  We first prove that $\eta^T\in C([0,1]\times [0,T])$. This assertion follows
  from the fact that $\bar{S}_{0T}\subseteq C([0,1]\times [0,T])$, which is proved
  as follows: First, for any $p>1$ by integrating the equation
$$
  \frac{d}{dt}\int_0^1|\eta(r,t)|^pr^2dr
  =p\int_0^1|\eta(r,t)|^{p-1}{\rm sgn}\eta(r,t)\eta_t(r,t)r^2dr
$$
  and using a similar argument as before we see that for any $\eta\in\bar{S}_{0T}$
  such that $\eta(\cdot,0)=0$, we have
$$
  \sup_{0\leq t\leq T}\int_0^1|\eta(r,t)|^pr^2dr
  =\frac{C\delta^{p-1}}{\sqrt{\varepsilon}}
  \Big(\int_0^T\!\!\int_0^1|\eta(r,t)|^2r^2drdt\Big)^{\frac{1}{2}}.
$$
  Next, by integrating the equation
$$
  \frac{d}{dr}|\eta(r,t)|^2=2\eta(r,t)\eta_r(r,t)
$$
  and using the H\"{o}lder inequality for $p>3$ we see that for any $\eta\in
  \bar{S}_{0T}$ we have
$$
\begin{array}{rl}
  \displaystyle\sup_{0\leq r\leq 1\atop 0<t<T}|\eta(r,t)|^2
  \leq &\displaystyle 2\sup_{0\leq t\leq T}\int_0^1|\eta(r,t)|dr
  \cdot\sup_{0\leq r\leq 1\atop 0<t<T}|\eta_r(r,t)|
  \leq C\delta\sup_{0\leq t\leq T}\int_0^1|\eta(r,t)|
  r^{\frac{2}{p}}\cdot r^{-\frac{2}{p}}dr
\\
  \leq &\displaystyle C\delta\sup_{0\leq t\leq T}
  \Big(\int_0^1|\eta(r,t)|^pr^2dr\Big)^{\frac{1}{p}}
\end{array}
$$
  Combining the above two inequalities we get
$$
  \sup_{0\leq r\leq 1\atop 0<t<T}|\eta(r,t)|\leq
  \frac{C\delta^{1-\frac{1}{2p}}}{\sqrt[p]{\varepsilon}}
  \Big(\int_0^T\!\!\!\int_0^1|\eta(r,t)|^2r^2drdt\Big)^{\frac{1}{4p}}.
\eqno{(6.6)}
$$
  Using this inequality to $\eta_n-\eta$ for any $\eta\in\bar{S}_{0T}$ and a
  corresponding sequence $\eta_n\in S_{0T}$ ($n=1,2,\cdots$) such that
  $\displaystyle\lim_{n\to\infty}\int_0^T\!\!\!\int_0^1|\eta_n-\eta|^2r^2drdt=0$,
  we see that $\eta_n$ converges to $\eta$ uniformly in $[0,1]\times [0,T]$, so
  that $\eta\in C([0,1]\times [0,T])$. This proves that $\bar{S}_{0T}\subseteq
  C([0,1]\times [0,T])$ and, consequently, $\eta^T\in C([0,1]\times [0,T])$.

  It follows that the problem (2.7) with $\eta$ replaced by $\eta^T$ has a unique
  classical solution for $0\leq t\leq T$, which we denotes as $(p^T,z^T)$. Since
  $\eta^T,p^T\in C([0,1]\times [0,T])$ and $z^T\in C[0,T]$, by using the standard
  $L^p$-theory for parabolic equations we see that when $(\eta,p,z)$ is replaced
  by $(\eta^T,p^T,z^T)$, the problem (2.8) has a unique strong solution in $[0,1]
  \times [0,T]$ which we denote as $\tilde{\eta}^T$. Clearly, $\tilde{\eta}^T\in
  W^{2,1}_p((0,1)\times (0,T))$ for any $1<p<\infty$. Using a limit argument we
  can easily show that $\tilde{\eta}^T=\eta^T$. Indeed, by choosing a sequence
  $\eta_n\in S_{0T}$ ($n=1,2,\cdots$) such that $\displaystyle\lim_{n\to\infty}
  \int_0^T\!\!\!\int_0^1|\eta_n-\eta^T|^2r^2drdt=0$ (which implies, by (6.6), that
  also $\displaystyle\lim_{n\to\infty}\sup_{0\leq r\leq 1\atop 0<t<T}|\eta_n-\eta^T|
  =0$) and defining $\tilde{\eta}_n=\mathscr{F}_T(\eta_n)$, we see that
  $\tilde{\eta}_n\to\tilde{\eta}^T$ uniformly in $[0,1]\times [0,T]$. Since also
  $\tilde{\eta}_n=\mathscr{F}_T(\eta_n)\to\bar{\mathscr{F}}_T(\eta^T)=\eta^T$
  strongly in $\bar{Y}_T$, by uniqueness of the limit we obtain the desired
  assertion. It follows that $(\eta^T,p^T,z^T)$ is a solution of the problem
  (2.6) in $[0,1]\times [0,T]$. Using this fact and the bootstrap argument we
  immediately obtain the assertion that $\eta^T_r\in C([0,1]\times [0,T])$.
  Hence $\eta^T\in S_{0T}$. The assertion that $\eta^T$ is a fixed point of
  $\mathscr{F}_T$ now follows from the fact that $\tilde{\eta}^T=\eta^T$.

  Since $(\eta^T,p^T,z^T)$ is a classical solution of the problem (2.6) in the
  time interval $[0,T]$, by uniqueness of the solution of this problem, we
  conclude that for any $0<T_1<T_2$ we have $\eta^{T_1}(r,t)=\eta^{T_2}(r,t)$
  for $(r,t)\in [0,1]\times [0,T_1]$. It follows that the following definition
  of the function $\eta$ in $[0,1]\times [0,\infty)$ makes sense:
$$
  \eta(r,t)=\eta^T(r,t) \quad \mbox{for any}\;\;(r,t)\in [0,1]\times [0,T]
  \;\; \mbox{and}\;\; T>0.
$$
  Moreover, by letting $T\to\infty$ in the relations $\eta^T\in S_{0T}$ and
  $\mathscr{F}_T(\eta^T)=\eta^T$, we see that $\eta\in S_0\subseteq S$ and
  it is a fixed point of the map $\mathscr{F}:S\to S$.

  Having proved that $\mathscr{F}$ has a fixed point in $S$, the desired
  assertions in Theorem 1.1 then immediately follow from Lemma 4.3 and
  Lemma 5.1. This completes the proof of Theorem 1.1. $\quad\Box$
\medskip


\end{document}